\documentclass[12pt,a4paper,reqno,oneside]{amsart}

\usepackage{t1enc}
\usepackage{amssymb}
\usepackage[mathscr]{euscript}

\usepackage{color}
\usepackage{float}

\usepackage{graphicx}
\usepackage{url}
\usepackage{caption}
\usepackage{subcaption}
\usepackage{dsfont}
\usepackage{gensymb}
\usepackage[toc,page]{appendix}

\usepackage{multicol}

\usepackage[capitalise]{cleveref}

\usepackage{lineno}

\usepackage[foot]{amsaddr}

\numberwithin{equation}{section}

\newtheorem{defi}{Definition}[section]
\newtheorem{thm}[defi]{Theorem}
\newtheorem{lemm}[defi]{Lemma}
\newtheorem{rem}[defi]{Remark}

\newtheorem{prop}[defi]{Proposition}

\newcommand{\clos}{\operatorname{cl}}
\newcommand{\convh}{\operatorname{conv}}

\newcommand{\B}{\mathcal{B}}
\newcommand{\EE}{\mathcal{E}}
\newcommand{\U}{\mathcal{U}}

\newcommand{\PP}{\mathbb{P}}

\newcommand{\E}{\mathbb{E}}

\newcommand{\N}{\mathbb{N}}
\newcommand{\R}{\mathbb{R}}

\newcommand{\TT}{\mathcal{T}}

\newcommand{\NN}{\mathcal{N}}
\newcommand{\FF}{\mathcal{F}}

\newcommand{\1}{\mathds{1}}

\begin{document}

\graphicspath{{plots/}}

\title[Environmental Contours for Non-Stationary Processes]
{Convex Environmental Contours for Non-Stationary Processes}

\author[Å. H. Sande]{Åsmund Hausken Sande$^*$}
\address{$^*$ Department of Mathematics, University of Oslo, Moltke Moes vei 35, P.O. Box 1053 Blindern, 0316 Oslo, Norway.}
\email{aasmunhs@math.uio.no}

\begin{abstract}
Environmental contours are tools frequently used in the early design of marine structures. They provide a description of critical design conditions and serve as a means for simplifying expensive long-term response calculations. Here, we consider convex contours based on the assumption of convex failure sets. We provide a rigorous foundation for the existence of such contours when the underlying environmental factors are modelled by a general, possibly non-stationary, process. This constitutes a generalisation of existing theory and is done to properly account for empirically observed increases in extreme sea-states.

Two definitions are proposed, based respectively on averages or quantiles of exceedence times, along with minimal conditions on the environmental processes to guarantee existence. In order to illustrate these methods we give two examples, including an empirical study containing a method for constructing contours based on the presented theory.
\end{abstract}

\maketitle



\section{Introduction}

Environmental contours are tools frequently used in the early design of marine structures. They provide a description of critical environmental conditions that may serve as a basis for structural design. Furthermore, they may be utilised to reduce the number of computationally expensive response calculations needed for the analysis of reliability. As a result, environmental contours are applied to analyse a wide variety of marine structures \cite{baarholm2010combining,fontaine2013reliability,giske2018long,vanemTrend}, and several methods for the construction of contours is mentioned in the \textit{recommended practices - environmental conditions and environmental loads} document by DNV (Det Norske Veritas) \cite{veritas2000environmental}.

A large variety of methods for constructing environmental contours exist, for a summary and comparison of different techniques we refer to \cite{contcompare} and \cite{contsummary}. A common thread through many of these methods is the modelling of environmental conditions, $V$, as a piecewise constant process with independent and identically distributed values in each interval. The length of this interval, $\Delta t$, varies depending on the application, but $\Delta t=3$ hours is a common choice. {The values of $V$ usually represent some summary statistics of the wave elevation, or other relevant environmental factors, over the period. $V$ then represents the long-term variations of the environmental conditions. The short-term variations, i.e.\ the variation of instantaneous conditions within the period, is usually ignored in the construction of these contours. }

These contours are also constructed to satisfy certain exceedence probabilities. These properties can usually be formulated by requiring that the probability of $V$ hitting a failure set $\FF$, not intersecting with the contour $\partial\B$, has at most a given probability $p_e$ of occurring in each interval. Due to the independence between different intervals, this assumption also implies restrictions on the time to failure $\tau_\FF=\inf\{t:V_t\in\FF\}$. Specifically, we have that the return period $\E[\tau_\FF]$ is bounded from below by $\Delta t/p_e$. Note that these exceedence properties are effectively defined under the assumption that failure depends only on $V$, thereby ignoring the short-term variation of the response.

Arguably, the most popular construction of environmental contours is the \textit{inverse first order reliability method} (IFORM) developed in \cite{IFORM,IFORM2008}. This method first establishes a sea-state distribution which induces a Rosenblatt transformation \cite{rosenblatt} of the density. All failure sets are then implicitly assumed to be convex in the transformed space, implying that the contour in the original space can be constructed by applying the inverse Rosenblatt transformation to a sphere.

As an alternative, in \cite{firstaltcontour,altcontour}, they develop a definition of convex environmental contours by assuming the failure set to be convex in the original space. This approach has several advantages, such as the easy inclusion of omission factors and a more amenable interpretation of the convexity assumption compared to IFORM. Several improvements and possible modifications to this method have been made in the literature. In \cite{dahl2018buffered} the concept of buffered contours is introduced and \cite{convcont} considers omission factors and convexity. Several different ways of constructing these contours are also discussed in e.g.\ \cite{voronoi,convcont,vanem20193,mackay2023model}.

Once a contour with the desired exceedence properties has been constructed, it can be applied to reliability analysis in several ways. Usually, response simulations are carried out for conditions along the contour over a period of length $\Delta t$. The point along the contour providing the highest extreme short-term response is then chosen as the design point. Often, the design point is chosen to correspond with the highest median of the response distribution. In e.g.\ \cite{giske2018long,sagrilo2011long}, an importance sampling procedure, centred around the design point, is discussed in the context of extreme long-term response computation. Additionally, in e.g. \cite{giske2018long,IFORM2008}, they consider quantiles of the response distribution at the design point as estimates of the characteristic response. These methods allow the estimation of response with only a limited number of computationally expensive response calculations.

However, these methods still rely on a stationary model. This causes issues when taken together with the evidence of increasingly extreme sea-states, as detailed in \cite{Kushnir,vanemTrend,vanemBayes}, which would imply a significant non-stationarity in significant wave heights.

Several articles such as \cite{vanemTrend} and \cite{firstaltcontour} adjust the sea-state distributions to correct for this increase. However, the models are still stationary, which keeps them from fully representing the changing behaviour of the environmental processes involved. A closer view on the differences between such strategies and the methods to be presented in this paper will be given in \cref{sec:empex}.

It is also worth mentioning the works of e.g.\ \cite{huseby2023AR1,Leira,mackay2023model,mackay2021effect,vanem2023analysing} which consider stationary processes with a varying degree of autodependence. These articles allow for more general behaviours of the underlying environmental processes, but do not address the issue of long-term trends.

The goal of this article is to present a mathematically rigorous framework for environmental contours with a broad class of possible models for the underlying environmental factors, including non-stationary ones. {While discrete models are the primary focus in terms of applications, we will also include the possibility of continuous-time models.} In this regard, we will give minimal conditions for relevant functions to be well defined in addition to existence of the contours themselves. 

Results presented in this article will be a generalisation of the theory discussed in \cite{convcont,firstaltcontour,altcontour}, which are based on the assumption of failure sets being convex in the original parameter space. As such, we will, in \Cref{sec:classiccontour}, give a brief overview of the main results and definitions from these papers, for so to generalise the setting in  \Cref{sec:dyncontour}. We will here propose two different ways of defining convex environmental contours based on either averages or quantiles of exceedence times. A brief discussion on how these definitions are connected to response analysis is given in \Cref{sec:response}. This is followed, in \Cref{sec:existence}, by a mathematically rigorous treatment of which minimal conditions are required of the model for $V$, in order to ensure the existence of contours. In \Cref{sec:thex} and \Cref{sec:empex} we present two examples of applications of the theory, which highlight some of the differences between classical approaches and the more flexible methods allowed by the theory presented in this article. As a part of the final example we also present a method for computing these contours in practice. 
 
\section{Convex Contours} \label{sec:classiccontour}

A \textit{convex environmental contour} is the boundary of a compact convex set $ \B \subset\R^d $, denoted $ \partial\B $, defined with respect to a $d$-dimensional environmental process $V$. For example, this process is often taken to be the pair $ V=(P,H) $ for $ d=2 $ where $ P $ is the zero-upcrossing wave period and $ H $ the significant wave height of a particular location of interest.

We will in this section follow the construction described in e.g.\  \cite{altcontour} and \cite{convcont} and assume that the distribution of $ V_t $ is constant and absolutely continuous with respect to the Lebesgue measure on $ \R^d $. Furthermore, the process is assumed be path-wise constant over periods of a set length of $ \Delta t $, and most importantly we make the assumption that values of $ V $ are independent between these different periods. One could equivalently consider $V_t=W_{\lfloor t/\Delta t\rfloor}$ where $\lfloor \cdot \rfloor$ denotes the floor function and with $\{W_n\}_{n=0}^\infty$ defined as a sequence of independent and identically distributed (i.i.d.) random variables with a distribution equal to that of $V$. As such we will refer to this type of model as an \textit{i.i.d. model} throughout this article.

{
For every possible structural design we consider a limit-state function $g$, also referred to as the performance function. This is assumed to depend only on $V$, thereby ignoring the variance of the structural response conditional on $V$. The function $g$ is defined such that the region $\FF$, where $g(V)\leq 0$, represents conditions the structure cannot safely handle. We therefore refer to $\FF$ as the \textit{failure set}. Environmental contours then aim to apply to any design satisfying $\FF\cap\B\subset\partial\B$. As such we consider an unknown performance function, and therefore an unknown failure set $\FF$. In order to handle such an unknown $ \FF $ we further assume that $ \FF $ belongs to $ \EE (\B) $, the class of all convex sets such that $ \FF\cap\B\subseteq\partial\B $. Based on this we may define the \textit{exceedence probability} by}
\begin{equation}\label{eq:exprob}
P_e(\B,\EE )=\sup_{\FF\in\EE (\B)}\PP(V\in\FF).
\end{equation}
and impose the constraint of 
\begin{equation}\label{eq:restriction0}
P_e(\B,\EE ) \leq p_e,
\end{equation}
where $ p_e $ is some given target exceedence probability. 

When dealing with convexity we will need the concept of hyperplanes. We will denote by $ \langle  \cdot , \cdot \rangle $ the canonical inner product on $ \R^d $ and by $ \|\cdot\| $ the euclidean norm, with this we also define the unit sphere on $ \R^d $ by $ S^{d-1}   = \{v\in\R^d:\|v\|=1\} $. The hyperplane indexed by the threshold $ c\in\R $ and the unit vector $ u\in S^{d-1} $ is then defined as
\begin{equation}
\Pi(u,c)   = \{v\in\R^d:\langle u,v\rangle  =  c \}.
\end{equation}

We further define the half-spaces
\begin{align}
	\begin{split}\label{def:halfspace}
		\Pi^-(u,c) &=  \{v\in\R^d:\langle u,v\rangle  \leq  c\},\\
		\Pi^+(u,c) &=  \{v\in\R^d:\langle u,v\rangle  \geq  c\},
	\end{split}
\end{align}
which allow us to present an important well-known result about separating hyperplanes. For a proof of this result, as well as \eqref{eq:BrepbyB}, we refer to \cite{rockafellar1997convex}.
\begin{prop}\label{prop:sephyper}
	For any two convex sets $ \B $ and $ \FF $ in $ \R^d $ such that $ \B\cap\FF\subseteq\partial\B $ there exists some $ u\in S^{d-1},\,c\in\R $ such that  $ \B\subseteq\Pi^-(u,c) $ and $ \FF\subseteq\Pi^+(u,c) $.
\end{prop}

This result implies that we can separate a convex set $ \B $ and any $ \FF\in\EE(\B) $. In particular, we can reconstruct any compact and convex $ \B $ as the intersection of all these half-spaces. If we define
\begin{equation}\label{def:Bu}
B(\B,u)=\sup\{\langle u,v\rangle: v\in\B\},
\end{equation}
we get that $ \B $ can be represented as
\begin{equation}\label{eq:BrepbyB}
	\B=\bigcap_{u\in S^{d-1}}\Pi^-(u,B(\B,u)).
\end{equation}

It is further shown in \cite{altcontour}, for $d=2$, that these hyperplanes also serve as maximal elements of $\EE(\B)$ for computing the exceedence probability. Specifically, we have
\begin{equation}\label{eq:exprobhalf}
P_e(\B,\EE )=\sup_{u\in S^{d-1}}\PP(V\in\Pi^+(u,B(\B,u))).
\end{equation}
Based on \eqref{eq:exprobhalf}, if we have an environmental contour $ \partial\B $ such that $ \B $ is convex and compact, we will call $ \partial\B $  a \textit{valid contour in the exceedence probability sense} if $ P_e(\B,\EE) \leq p_e$, and a \textit{proper contour in the exceedence probability sense} if, for all $ u\in S^{d-1} $, $ \PP(V\in \Pi^+(u,B(\B,u))) = p_e $. The goal is then to construct the smallest convex and compact set with a valid or, ideally, proper contour.

If we then define
\begin{equation}\label{def:Cu}
C_e(u)=\inf\{C:\PP( \langle u,V\rangle  >C)=p_e\},
\end{equation}
and if $ B(\B,u)\geq C_e(u) $ for all $ u\in S^{d-1} $ we get $ P_e(\B,\EE ) \leq p_e $, making $ \partial\B $ a valid contour. In particular if there exists any convex $ \B $ such that $ B(\B,\cdot)=C_e(\cdot) $ then $\partial\B$ is the unique proper contour. Lastly, if this is the case, then 
\begin{equation}\label{eq:CeqB_def}
\B=\bigcap_{u\in S^{d-1}}\Pi^-(u,C_e(u)),
\end{equation}
gives an explicit construction of this optimal proper contour.

In the following sections we will need the concept of the first hitting time of a set $ \FF\subseteq\R^d $, defined by
\begin{equation}
	\tau_\FF = \inf\{t:V_t\in\FF\}	
\end{equation}

For our i.i.d.\ model we can easily associate our target exceedence probability with a target return period. We can note that by our assumption of independence between the $ W_n $'s that the exceedence time, $\tau_\FF$, is geometrically distributed. Furthermore, for any valid $\partial\B$ and $\FF\in\EE(\B)$, the mean of $\tau_\FF$ is at least $ \Delta t/p_e $. This implies that when we are ensuring that $ P_e(\B,\EE )\leq p_e $ we are equivalently ensuring that $ \E[\tau_\FF]\geq t_r $ for some target return period $ t_r = \Delta t/p_e$. Similar arguments would allow us to compare $ P_e(\B,\EE ) $ with quantiles of the distribution of $ \tau_\FF $. Both the mean and quantiles of $ \tau_\FF $ are more amenable to generalisation than exceedence probabilities and will be used in the upcoming sections.

\section{Environmental Contours for General Processes} \label{sec:dyncontour}

We now aim to extend the concepts introduced in the previous section to a more general context. We no longer assume $V$ to be stationary and instead consider it to be a progressively measurable process taking values in $\R^d$. We also need the process to satisfy certain regularity conditions in order to ensure that \eqref{eq:phidef}, which will be introduced later, is measurable. For this purpose one may assume, for instance, càdlàg paths. {Usually, discrete models are used in order to facilitate response analysis. Fortunately, a discrete model for $V$ is sufficient to ensure the measurability of \eqref{eq:phidef}.}

We will also still consider an unknown failure set $ \FF\in\EE (\B) $ where $ \EE (\B) $ is the collection of all convex sets $ \FF $ such that $ \FF\cap\B\subseteq\partial\B $. This will similarly allow the use of half-spaces to control the exceedence time. 

Since we no longer assume a stationary distribution, we will introduce two ways of replacing \eqref{eq:exprob}. A common substitute for the exceedence probability, used explicitly in works such as \cite{huseby2023AR1,Leira,mackay2023model,mackay2021effect,vanem2023analysing}, is to use the average failure time, commonly referred to as the \textit{return period}. As such, we start by defining the \textit{return period of $ \B $} by
\begin{equation}\label{eq:returnperiod}
	T_r(\B) = \inf_{\FF\in\EE (\B)}\E[\tau_\FF].
\end{equation}
\begin{rem}
	Since $ V $ is now possibly non-stationary the concept of a long-term average return period is no longer meaningful. However, for the sake of consistency, we shall still refer to these average exceedence times as return periods.
\end{rem}
In some cases there may be yearly trends present which, if persisting indefinitely, may cause the process to drift over time. In such cases there could be a positive probability that the process never exits certain sets. This contrasts with the stationary case where every set of positive measure (w.r.t. the law of $V$) is eventually hit. While the exceedence time might have a positive probability of not occurring, thereby making the return period infinite, there could still be a high chance of it occurring in finite time. In order to account for such behaviour we want a more flexible version of \eqref{eq:returnperiod}. As such we define the \textit{survival probability of $ \B $} by
\begin{equation}
Q_s(\B) = \inf_{\FF\in\EE (\B)}\PP(\tau_\FF > t_s),
\end{equation}
for a given \textit{survival time} $ t_s>0 $.

Like with \eqref{eq:restriction0}, we will construct our contour based on these two definitions. In our case we will consider two separate possible restrictions, the first one is based on return periods with 
\begin{equation}\label{eq:restriction1}
T_r(\B) \geq t_r 
\end{equation}
for some target return period $ t_r >0 $. If this holds then for any $ \FF\in\EE (\B) $ it will take on average at least a time of $ t_r $ to enter $ \FF $. Note that under the constrains of an i.i.d\ model, \eqref{eq:restriction1} is equivalent to \eqref{eq:restriction0} for $ t_r = \Delta t/p_e$.

The alternative restriction corresponding to the survival probability is defined as
\begin{equation}\label{eq:restriction2}
	Q_s(\B) \geq q_s	
\end{equation}
for a given minimal survival probability $1 > q_s > 0$. If this condition holds, then for any $ \FF\in\EE(\B) $, it is guaranteed that with a probability of at least $ q_s $, the process $V$ will take at least $ t_s $ amount of time before hitting $ \FF $. Note that under the constraints of an i.i.d\ model, $ \tau $ has a geometric distribution which means \eqref{eq:restriction2} is equivalent to \eqref{eq:restriction0} if e.g.\ $ t_s = \Delta t/p_e$ and $ q_s=(1-p_e)^{1/p_e} $. In particular, for low exceedence probabilities, we have that $ q_s\approx 1/e \approx 37\% $.

{
\begin{rem}
    While contours are usually formulated using return periods, there are several benefits to considering contours based on \eqref{eq:restriction2}. As mentioned, long-term trends can lead to heavy tails of $\tau_\FF$ which can inflate $\E[\tau_\FF]$, or even make it infinite. There are also two technical benefits in that only the path of $V$ up to time $t_s$ needs to be considered. If these paths need to be simulated, as in \Cref{sec:empex} or \cite{vanem2023analysing}, then simulating the average exceedence time could require paths of arbitrary length, introducing certain numerical challenges. Lastly, given the difficulty in forecasting distant future trends, avoiding the specification of these trends for timepoints beyond $t_s$ is advantageous.
\end{rem}
}

With these possible constraints established we now have two new ways of defining our environmental contours, either by \eqref{eq:restriction1} or \eqref{eq:restriction2}. As noted, both of these serve as generalizations of the restriction of \eqref{eq:restriction0}. Analogously to the previous section we will refer to a contour $ \partial\B $ such that $ \B $ is convex and compact as \textit{valid in the return period sense} if $ T_r(\B) \geq t_r  $ and \textit{valid in the quantile sense} if $ Q_s(\B) \geq q_s $. Likewise, we call $ \partial\B $ \textit{proper in the return period sense} if $ \E[\tau_{\Pi^+(u,B(\B,u))}]=t_r $ for all $ u\in S^{d-1} $ and \textit{proper in the quantile sense} if $ \PP(\tau_{\Pi^+(u,B(\B,u))} > t_s)=q_s $ for all $ u\in S^{d-1} $. In order to justify the definitions of proper contours we proceed analogously to \cite{altcontour} by the following result. 
\begin{prop}\label{prop:maxhalfspace}
	Let $\B\subset\R^d$ be a compact and convex set, we then have
	\begin{align*}
		Q_s(\B) &= \inf_{u\in S^{d-1} }\PP(\tau_{\Pi^+(u,B(\B,u))} > t_s),\\
		T_r(\B) &= \inf_{u\in S^{d-1}}\E[\tau_{\Pi^+(u,B(\B,u))}].
	\end{align*}

	\begin{proof}
		
		We first note that by \cref{prop:sephyper} we have for any $ \FF\in\EE (\B) $ that $ \FF\subseteq\Pi^+(u,c) $ and $ \B\subseteq \Pi^-(u,c)$ for some $ u\in S^{d-1},\,c\in\R$. This yields $ B(\B,u) \leq c$, which implies $ \Pi^+(u,c)\subseteq \Pi^+(u,B(\B,u)) $. Since $ \FF\subseteq\Pi^+(u,B(\B,u)) $ we then have $ \tau_\FF \geq \tau_{\Pi^+(u,B(\B,u))} $ which further implies that  
		\begin{align*}
			\PP(\tau_{\FF} > t_r) &\geq
			\PP(\tau_{\Pi^+(u,B(\B,u))} > t_r), \\ 
			\E[\tau_{\FF}]        &\geq
			\E[\tau_{\Pi^+(u,B(\B,u))}].
		\end{align*}
		
		This inequality yields
		\begin{align*}
			\inf_{\FF\in\EE (\B)}\PP(\tau_\FF \geq t_s)
			& \geq \inf_{u\in S^{d-1} }\PP(\tau_{\Pi^+(u,B(\B,u))} > t_r),\\
			\inf_{\FF\in\EE (\B)}\E[\tau_\FF]
			& \geq \inf_{u\in S^{d-1}}\E[\tau_{\Pi^+(u,B(\B,u))}].
		\end{align*}
		The result then follows by noting that $ \Pi^+(u,B(\B,u))\in\EE (\B) $ for any $u\in S^{d-1}$.
	\end{proof}
\end{prop}

\begin{rem}
    \Cref{prop:maxhalfspace} treats the tangent half-spaces,  $\Pi^+(\cdot,B(\B,\cdot))$, as maximal elements of $\EE(\B)$. However, it is also possible to interpret the tangent half-spaces as FORM approximations of failure sets for possible designs. This would replace our assumption of convex failure sets by the linear FORM approximation.
\end{rem}

With \Cref{prop:maxhalfspace}, we can introduce our analogues of $ C_e $ from \eqref{def:Cu} by
\begin{align}
		C_Q(u)&=\inf\{b\in\R:
		\PP\left(\tau_{\Pi^+(u,b)} > t_s\right) = q_s \} , \label{def:CQ}\\
		C_T(u)&=\inf\{b\in\R:
		\E\left[\tau_{\Pi^+(u,b)} \right] = t_r \}.\label{def:CT}
\end{align}
These functions allow us to define our contours by
\begin{align*}
    Q_s(\B)\geq p_s &\Leftrightarrow B(\B,u) \geq C_Q(u) \text{ for all }u\in S^{d-1},\\
	T_r(\B)\geq t_r &\Leftrightarrow B(\B,u) \geq C_T(u) \text{ for all }u\in S^{d-1},
\end{align*}
which will be proved in \Cref{prop:BgeqC_properties}. However, before we move to the theoretical considerations of this article we will discuss the more practical connection to response analysis.

\section{Connection with Response}\label{sec:response}

\subsection{Interpretation of Environmental Processes in Continuous Time}

The theory presented makes no restrictions on whether $V$ is modeled as discrete or not, which causes some minor complications in applying these contours to response analysis. 

In this regard it is worth mentioning that the use of continuous-time processes for the definition of contours has been previously considered in e.g.\ \cite{Leira}. Here, a continuous process is made comparable to an i.i.d.\ process by equating the outcrossing rate over a period of length $T$ for specific thresholds. This was achieved through the use of \textit{equivalent characteristic time scales}. This procedure ensured that the continuous and i.i.d.\ process produced the same contour for a given return period, but whenever the target return period was changed the resulting contours differed significantly.

For offshore engineering it is common to split the description of ocean waves into its long-term and short-term variability. The long-term variability often considers summary statistics such as significant wave height and zero-upcrossing period, these describe conditions over a certain time period $\Delta t$. For example, the significant wave height over this period would be the average height of the largest third of waves within this period. The short-term variation usually describes the variation of individual waves within those periods.

Generally, the environmental conditions along a contour represents the long-term conditions, i.e.\ summary statistics over a period of length $\Delta t$. Response analysis is the carried out by simulating the short-term variations over a period of that length, under the assumption of constant long-term conditions, see e.g.\ \cite{baarholm2010combining,giske2018long, sagrilo2011long}. In order to facilitate this type of response analysis it is important that we can interpret $V$ as suitable summary statistics of the sea-state over some period. There are several ways of ensuring this, most easily and commonly done by modelling $V$ as a discrete process. 

{If one wanted to model $V$ as a continuous-time process it is possible to consider $V_{t}$ to represent some collection of summary statistics over $[ t,t+\Delta t]$. This can be achieved by e.g.\ choosing a continuous-time model equating the distribution of $V_{t}$ to an estimated density of the relevant statistics over the given period. In so doing one must also ensure that the autodependence structure of $V$ is sufficiently accurate. An simple example of this approach is given in \Cref{sec:thex}.
}

\subsection{Applications to Response Analysis}

Environmental contours are often applied in the early concept evaluation phase, where the contours help identify possible critical design conditions. In terms of applications to response analysis we mainly have deterministic response, characteristic response estimates by way of design points, and importance sampling centred around these design points. We will here discuss the two former, as the extension of the latter to a nonstationary setting is still quite speculative.

\subsubsection{Deterministic Response}\label{sec:detresponse}

The most straightforward applications of contours is in the case of deterministic response. While usually unrealistic in practise, this case highlights the intuitive basis for the use of environmental contours.

Assume we have a deterministic response function $y$ and a target return period $t_r$ for the design of our structure. For any design we will then have a response capacity, $y_\text{cap}$ inducing a limit-state function, $g(v)=y_\text{cap}-y(v)$, and an associated failure set $\FF=\{v\in\R^d:g(v)<0\}$. The goal is then to find a minimal response capacity, $y_\text{cap}$, such that the mean time to failure, $\E[\tau_\FF]$, is at least $T$ years. This can easily be done by considering any valid contour in the return period sense, as defined in \Cref{sec:dyncontour} for a $T$-year return period. We can then apply an analogue to the \textit{inverse FORM method} described in \cite{IFORM}, which chooses 
$$y_\text{cap}=\max_{b\in\partial\B}y(b).$$
Assuming that $\text{argmax}_{b\in\B} \, y(b)\in\partial\B$ and that the resulting $\FF$ is convex, we then have by the definition of the contour that $\E[\tau_\FF]\geq T$ years. Since, in the simple deterministic case, the failure of the structure occurs at time $\tau_\FF$, we know that the average failure time of the structure shares the same bound.

\subsubsection{Quantiles of Design Point}
{The most common application of environmental contours in the case of stochastic response}, mentioned in e.g.\ \cite{giske2018long,IFORM2008}, is the following. First compute an environmental contour with a $T$-year return period. For environmental conditions along the contour, compute the distribution of the short-term maximal response. Designate the worst of these conditions as the design point, the conditions providing a response distribution with the highest median is often chosen. The response level with a $T$-year return period is then estimated by a quantile of the response distribution at the design point, quantiles ranging from $85\%$ to $95\%$ is recommended in \cite{NORSOK2017}.  This procedure, as described, can also be applied to a non-stationary setting.

{Note that by \Cref{sec:detresponse}, the design point estimate would correspond to an exact bound in the deterministic case. The idea is then to choose a suitably high quantile to correct for the stochastic nature of the response.}
The underlying assumption then becomes that a representative response value with a return period of $T$ years should occur along a contour with the same return period (paraphrasing \cite{Leira,contsummary}). By applying an analogous assumption to our non-stationary contours then this can shed some light on the choice between contours based on survival times or return periods. {If one is interested in finding a response level with a specified average exeedance time, then a contour based on that same average exceedence time should be used. Similarly, if quantiles of this exceedence time is of primary interest, then contours based on survival times should be considered.} 

It is important to note that the use of quantiles of the design point for the calculation of characteristic response levels are rough approximations. It is usually recommended to verify the choice of quantile by a full long-term response analysis if possible. Despite this, the method is highly efficient and requires very few response simulations to be carried out. 

{In the case where $V$ is modeled in continuous time it is still possible to identify a design point along the contour. Assuming that $V$ represents {the long-term conditions} over a period of length $\Delta t$, then response distributions can be established for conditions along the contour. However, since a full response analysis in continuous time is generally unfeasible, there is no practical way to verify the choice of quantile. In the case where a full response analysis is out of reach, continuous-time modelling is an option. However, this issue is a strong reason to focus on discrete models whenever a full response analysis is needed. }

\section{Existence of Contours}\label{sec:existence}
The design conditions along the contour are chosen due to their statistical properties, as such it is important to be able to clearly interpret and mathematically verify them. In \Cref{sec:dyncontour}, we allow for a very general class of models for $V$, however, not every model permits the existence of well-defined contours. Therefore, in this section we will provide rigorous mathematical justification for the existence of these contours. Firstly, we give minimal conditions for $C_Q$ and $C_T$ to be well defined. We will then show that the analogous representation of \Cref{eq:CeqB_def}, based on constructing contours by $B(\B,u)=C_e(u)$, still provides a unique proper contour in our generalised setting. Finally, in the case where no proper contours exist, we prove existence of valid contours.

In order to ensure our functions are well defined we will need the following definitions and results. For any $ u\in S^{d-1} $, $b\in\R$ and $t\in\R $ with $ t\geq0 $, we denote the cumulative distribution function of $\tau_{\Pi^+(u,b)}$ by
\begin{equation}
	F^u_b(t)=\PP\left(\tau_{\Pi^+(u,b)} \leq t \right),
\end{equation}
and the average of $\tau_{\Pi^+(u,b)}$ by
\begin{equation}
	\TT_u(b)=\E\left[\tau_{\Pi^+(u,b)}\right].
\end{equation}
Finally we define
\begin{equation}\label{eq:phidef}
	\phi^u_t = \sup_{s\in [0,t]}\langle V_s,u\rangle.
\end{equation}

Throughout this article we will usually assume that for any $t>0$, $ b\in\R\cup\{\pm\infty\} $  we have
\begin{equation}\label{axiom:contmax}
	\PP\left( \phi^u_t = b\right) = 0.
\end{equation}
This assumption will serve as a minimal condition for our contours and other concepts to be definable. For most models, \eqref{axiom:contmax} will follow as a consequence of $ \phi^u_t $ admitting a continuous density for every $ u\in S^{d-1} $ and $ t>0 $. For a discrete model of $V$, it is sufficient that $V_{t}$ admits a density for all $t$. 

To see the connection between these definitions we have the following result.

\begin{lemm}\label{lem:phieq}
	Let $ u\in S^{d-1} $ and $ t>0 $, we then have
	$$\TT_u(b)=\int_{0}^{\infty}\left(1-F_b^u(t)\right)dt.$$
	Furthermore, if \eqref{axiom:contmax} holds, then
	$$ F^u_b(t)  =  \PP\left(\phi^u_t \geq b\right).$$
	\begin{proof}
		The first equality is the standard tail probability expectation formula, as such we omit the proof.
		
		As for the second equality, if $\tau_{\Pi^+(u,b)} \leq t$ then there is some point $s\leq t$ such that $V_s\in\Pi^+(u,b)$, or equivalently $\langle V_s,u\rangle \geq b$, which implies $\phi^u_t \geq b$. Similarly, if $\tau_{\Pi^+(u,b)} > t$ then no such point exists and consequently $\phi^u_t \leq b$, this conversely states that $\phi^u_t > b$ implies $\tau_{\Pi^+(u,b)} \leq t$.
		
		Applying these implications and \eqref{axiom:contmax} we get
		\begin{equation*}
			\PP\left(\tau_{\Pi^+(u,b)} \leq t \right)
			\leq  \PP\left(\phi^u_t \geq b \right)
			=  \PP\left(\phi^u_t > b \right)
			\leq \PP\left(\tau_{\Pi^+(u,b)} \leq t \right),
		\end{equation*}
		which proves the second equality.
	\end{proof}
\end{lemm}

With this lemma we can prove the following results which guarantee that \eqref{def:CQ} and $\eqref{def:CT}$, i.e.\ $C_Q$ and $C_T$, are well defined.

\begin{thm}\label{thm:Fbtcont}
	Under the assumption of \eqref{axiom:contmax} we have for any $ t > 0 $, $u\in S^{d-1}$ that $ b\mapsto F^u_b(t) $ is monotone non-increasing and continuous with $ \lim_{b\to\infty}F^u_b(t)=0 $ and $ \lim_{b\to-\infty}F^u_b(t)=1 $.
	\begin{proof}
		We will start by showing that $ F $ is monotone non-increasing in $ b $. If $ b\leq b' $ then $ \Pi^+(u,b') \subseteq \Pi^+(u,b) $ which implies $ \tau_{\Pi^+(u,b)} \leq  \tau_{\Pi^+(u,b')} $, finally yielding that $ F^u_{b'}(t)\leq F^u_{b}(t)$.
		
		For the continuity we start by showing left-continuity of $ b\mapsto F^u_b(t) $. Whenever $ b_n\to b,\, b_n < b_{n+1} $ we have by set-continuity of measures and \Cref{lem:phieq} that
		\begin{align*}
			\lim_{n\to\infty}F^u_{b_n}(t)-F^u_b(t)	
			=& \lim_{n\to\infty} \PP\left(\tau_{\Pi^+(u,b_n)} \leq t \right)
			-\PP\left(\tau_{\Pi^+(u,b)} \leq t \right) \\
			=& \PP\left(\bigcap_n\{\tau_{\Pi^+(u,b_n)} \leq t\} \right)
			-\PP\left(\tau_{\Pi^+(u,b)} \leq t \right) \\
			=& \PP\left(\tau_{\Pi^+(u,b_n)} \leq t \,\text{ for all }n\right)
			-\PP\left(\tau_{\Pi^+(u,b)} \leq t \right) \\
			=& \PP\left(\phi^u_t \geq b_n \,\text{ for all }n\right)
			-\PP\left(\phi^u_t \geq b \right).	
		\end{align*}

		 We see that $ \phi^u_t \geq b_n $ for all $ n $, is equivalent to $ \phi^u_t \geq b $. This implies that the limit equals $\PP\left(\phi^u_t \geq b \right)-\PP\left(\phi^u_t \geq b \right)=0$, thus proving the left-continuity of $ b\mapsto F^u_b(t) $.

		As for right-continuity we consider $ b_n\to b,\, b_n > b_{n+1} $ and get
		\begin{align*}
			F^u_b(t)-\lim_{n\to\infty}F^u_{b_n}(t)
			&= \PP\left(\tau_{\Pi^+(u,b)} \leq t \right)
			-\lim_{n\to\infty} \PP\left(\tau_{\Pi^+(u,b_n)} \leq t \right)
			\\
			&= \PP\left(\tau_{\Pi^+(u,b)} \leq t \right)
			- \PP\left(\bigcup_n\{\tau_{\Pi^+(u,b_n)} \leq t\} \right)
			\\
			&= \PP\left(\phi^u_t \geq b \right)
			- \PP\left(\bigcup_n\{\phi^u_t \geq b_n\} \right)
			\\
			&= \PP\left(\phi^u_t \geq b \right)
			- 1 + \PP\left(\bigcap_n\{\phi^u_t < b_n\} \right)
			\\
			&= \PP\left(\phi^u_t > b \right)
			 + \PP\left(\phi^u_t < b_n \,\text{ for all }n \right) - 1
			\\.
		\end{align*}
		Similarly to before we note that $ \phi^u_t < b_n $ for all $ n $ is equivalent to $ \phi^u_t \leq b $. As a consequence the limit equals $\PP\left(\phi^u_t \in\R \right)-1=0$, which implies that $ b\mapsto F^u_b(t) $ is right-continuous and therefore fully continuous.
		
		Considering $ \lim_{b\to\infty}F^u_b(t) $ we may, for any sequence $ b_n\to \infty,\, b_n<b_{n+1} $, compute
		\begin{align*}
			\lim_{n\to\infty} F^u_{b_n}(t)  
			&= \lim_{n\to\infty}\PP(\tau_{\Pi^+(u,b_n)}\leq t)\\
			&= \lim_{n\to\infty}\PP(\phi^u_t \geq b_n)\\
			&=  \PP(\phi^u_t \geq b_n \,\text{ for all }n)\\
			&=  \PP(\phi^u_t = \infty)\\
			&= 0.
		\end{align*}

		And lastly,  for $ \lim_{b\to-\infty}F^u_b(t) $ we get, for any sequence $ b_n\to -\infty,\, b_n>b_{n+1} $, that
		\begin{align*}
		\lim_{n\to\infty} F^u_{b_n}(t)
		&= \lim_{n\to\infty}\PP\left(\tau_{\Pi^+(u,b_n)}\leq t\right)\\  
		&= 1-\lim_{n\to\infty}\PP\left(\tau_{\Pi^+(u,b_n)}> t\right)\\  
		&= 1-\lim_{n\to\infty}\PP\left(\phi^u_t\leq b_n\right)\\
		&=  1- \PP\left(\phi^u_t \leq b_n \,\text{ for all }n\right)\\
		&=  1- \PP\left(\phi^u_t = -\infty\right)\\
		&=  1,		
		\end{align*}
		which completes the proof.
	\end{proof}
\end{thm}

\begin{rem}
	
	Continuity of $ b\mapsto F^u_b(t) $ and $ \lim_{b\to\infty} F^u_{b}(t)=0 $ for all $u\in S^{d-1}$ is equivalent to the assumption of \eqref{axiom:contmax}, making it a necessary and sufficient condition for \Cref{thm:Fbtcont}.
\end{rem}

This theorem implies that $ b\mapsto 1-F^u_b(t_s)=\PP(\tau_{\Pi^+(u,b)}>t_s) $ spans $ (0,1) $ which implies that $ C_Q $ is well defined for any $ q_s\in (0,1) $. To ensure that $ C_T $ is also well defined, we have the following result.

\begin{prop}
	Assume that \eqref{axiom:contmax} holds and that for any $u\in S^{d-1}$ there is some $ b_u^*\in\R\cup\{\infty\} $ such that $ \TT_u(b)<\infty $ for all $ b\in (-\infty,b_u^*) $ with $ \TT_u(b)=\infty $ for all $b\geq b_u^* $. We then have that $ \TT_u(\cdot) $ is continuous and monotone non-decreasing on $ (-\infty,b_u^*) $ with $ \lim_{b\to -\infty}\TT_u(b)=0 $ and $ \lim_{b\to b_u^*}\TT_u(b)=\infty $.
	\begin{proof}
		We start with monotonicity. If $ b\leq b' $ then $ \Pi^+(u,b') \subseteq \Pi^+(u,b) $  which means $ \tau_{\Pi^+(u,b)} \leq \tau_{\Pi^+(u,b')}$, implying $ \TT_u(b)\leq \TT_u(b') $.
		
		As for continuity, if $ b_n \to b\in (-\infty,b_u^*) $ there exists some $ \epsilon>0 $ and some $ N\in\N $ such that $ b_n < b+\epsilon $ for all $ n > N $ with $ b+\epsilon \in (-\infty,b_u^*) $. This means that $ 1-F^u_{b_n}(t) \leq 1-F^u_{b+\epsilon}(t)$ for all $ n>N $ and $ \TT_u(b+\epsilon)=\int_{0}^{\infty}(1-F^u_{b+\epsilon}(t))dt<\infty $. We then get by continuity of $ b\mapsto F^u_b(t) $ and the dominated convergence theorem that
		\begin{align*}
			\lim_{n\to \infty}\TT_u(b_n)
			=&\lim_{n\to \infty}\int_{0}^{\infty}\left(1-F^u_{b_n}(t)\right)dt\\
			=&\int_{0}^{\infty}\left(1-\lim_{n\to \infty}F^u_{b_n}(t)\right)dt\\
			=&\int_{0}^{\infty}(1-F^u_{b}(t))dt\\
			=&\TT_u(b). 
		\end{align*}
		
		Similarly, for $ b_n \to -\infty $, we have some $ N $ and $b'\in (-\infty,b_u^*)$ such that $ b_n < b' $ for all $ n>N $. Since $ 1-F^u_{b_n}(t)\leq 1-F^u_{b'}(t) $ for $n>N$ we get by the dominated convergence theorem, along with $ \lim_{b\to -\infty}F^u_{b}(t)=1 $ , that
		\begin{align*}
		\lim_{b\to -\infty}\TT_u(b)
		=&\lim_{b\to -\infty}\int_{0}^{\infty}(1-F^u_{b})(t)dt\\
		=&\int_{0}^{\infty}\left(1-\lim_{b\to -\infty}F^u_{b}(t)\right)dt\\
		=&\int_{0}^{\infty}0\,dt\\
		=&\,0. 
		\end{align*}
		
		Lastly, for $ b_n \to b_u^*<\infty, $ we start with right limits and assume  $ b_n  \geq b_{n+1} $ for all $ n $ with $ b_n\to b_u^*<\infty $. Since $ b_n \geq b_u^*$ we have $ \TT_u(b_n) =\infty $ for all $ n $ which yields $ \TT_u(b_n)\to \infty =\TT_u(b)$.
		
		We then consider the left limit case, i.e.\  $\, b_n \leq b_{n+1} $ for all $ n $, $b_n\to b_u^*$. We then have that $ \{1-F^u_{b_n}\}_{n=1}^\infty $ is a monotone increasing sequence of non-negative functions, as such we get by the monotone convergence theorem that
		\begin{align*}
		\lim_{b\to b_u^*}\TT_u(b)
		=&\lim_{b\to b_u^*}\int_{0}^{\infty}\left(1-F^u_{b}(t)\right)dt\\
		=&\int_{0}^{\infty}\left(1-\lim_{b\to b_u^*}F^u_{b}(t)\right)dt\\
		=&\int_{0}^{\infty}(1-F^u_{b_u^*}(t))dt\\
		=&\,\TT_u\left(b_u^*\right)\\
		=&\,\infty.
		\end{align*} 
		The same computations would hold if $ b_u^*=\infty $ by considering $ F^u_{\infty}=0 $, which completes the proof.
	\end{proof}
\end{prop}

With this result we see that, under the given assumptions, $ \TT_u(\cdot) $ spans the whole of $ (0,\infty) $ so $ C_T $ is well defined for any $ t_r \in (0,\infty)$. 

With this, both our analogues of $ C_e $ from \eqref{def:Cu} are well defined. Similarly to \cite{altcontour}, we can use these functions to guarantee certain properties of our contours.

\begin{prop}\label{prop:BgeqC_properties}
	Fix some $ t_s \in (0,\infty),\,q_s\in (0,1) $ such that $ C_Q(u) $ is well defined for all $ u\in S^{d-1} $. We then have that $ Q_s(\B)\geq p_s  $ is equivalent to
	$$B(\B,u) \geq C_Q(u) \text{ for all }u\in S^{d-1}.$$
	Furthermore, if we fix some $t_r > 0$ such that $C_T$ is well defined, then $ T_r(\B)\geq t_r $ is equivalent to 
	$$B(\B,u) \geq C_T(u) \text{ for all }u\in S^{d-1}.$$
	\begin{proof}
		We start by assuming that $ B(\B,u) \geq C_Q(u) \text{ for all }u\in S^{d-1} $ and get
		$$ Q_s(\B)
		= \inf_{u\in S^{d-1}} \left\{\PP\left(\tau_{\Pi^+(u,B(\B,u))} > t_s\right)\right\} 
		\geq \inf_{u\in S^{d-1}}\left\{\PP\left(\tau_{\Pi^+(u,C_Q(u))} > t_s\right)\right\}
		 = q_s.  $$
		
		Conversely, if $ B(\B,u') < C_Q(u') \text{ for some }u'\in S^{d-1} $ then by the definition of $ C_Q $ we must have $\PP\left(\tau_{\Pi^+(u',B(u'))} > t_s\right) <  p_s $. This implies
		$$ Q_s(\B)
		= \inf_{u\in S^{d-1}} \left\{\PP\left(\tau_{\Pi^+(u,B(\B,u))} > t_s\right)\right\} 
		\leq \PP\left(\tau_{\Pi^+(u',B(\B,u'))} > t_s\right)
		< q_s.  $$
		An identical argument proves the statement about $ T_r $ which completes the proof.
	\end{proof}
\end{prop}

Recall that we can write 
\begin{equation}
\B=\bigcap_{u\in S^{d-1}}\Pi^-(u,B(\B,u)),
\end{equation}
so if there exists some $ \B $ with $ B(\B,u)=C(u) $ for all $ u\in S^{d-1} $, then we can immediately construct $ \B $ by 
\begin{equation}\label{eq:propconstr}
\B=\bigcap_{u\in S^{d-1}}\Pi^-(u,C(u)),
\end{equation}

for e.g.\ $C=C_e$, $C_T$ or $C_Q$. In \cite{altcontour} it is shown for an i.i.d.\ model of $ V $ that, under some conditions, $ B(\B,u)=C_e(u) $ is equivalent to $ \partial\B $ being proper in the exceedence probability sense. This implies that all proper contours are constructable in the same fashion as \eqref{eq:propconstr}. An analogous result also holds in our setting for $C_Q$ or $C_T$ under some light assumptions.

\begin{prop}\label{prop:propQ}
	For any $ u\in S^{d-1} $, $ t\in(0,\infty) $, define the set $ \U^u_{t} = \{b:F^u_b(t)\in(0,1)\} $ and assume that $ \phi^u_t $ admits a density, denoted by $ f^u_t $, satisfying $  f^u_t(b)>0  $ for all $ b\in\U^u_{t} $. Note that this last requirement is equivalent to the support of $ f^u_t $ being a connected interval.
	
	We have that $ b\mapsto F^u_b(t) $ is monotone decreasing on $ \U^u_{t} $. Furthermore, if there exists some proper contour $ \partial\B $ in the quantile sense, i.e.\ $ \PP(\tau_{\Pi^+(u,B(\B,u))} > t_s)=q_s $ for all $ u\in S^{d-1} $, then $ B(\B,\cdot)=C_Q $ and therefore
	$$\B=\bigcap_{u\in S^{d-1}}\Pi^-(u,C_Q(u)).$$

	\begin{proof}
		Since $ \phi^u_t $ admits a density we have that $ \PP\left( \phi^u_t = b\right) = 0 $ for all $ b\in\R $ which implies by \Cref{thm:Fbtcont} that $ b\mapsto F^u_b(t) $ is continuous and monotone non-increasing, in particular we note that $ \U^u_{t} $ is open for all $ u\in S^{d-1} $, $ t\in(0,\infty) $.
		
		We first aim to prove that $ b\mapsto F^u_b(t) $ is monotone decreasing on $ \U^u_{t} $. To see this we consider $ b\in \U^u_{t} $ and $ b'\in\R $ such that $ b < b' $. Since $ \U^u_{t} $ is open we can find an $ \epsilon > 0 $ such that $ (b,b+\epsilon)\subseteq\U^u_{t} $ and $ b+\epsilon < b' $. Additionally, we have
		$$ \PP\left(\phi_t^u \in (b,b+\epsilon)\right)
		= \int_{b}^{b+\epsilon}f^u_t(x)dx
		>0.$$ 
		Combining this with \cref{lem:phieq} then yields
		\begin{align*}
			F^u_b(t)-F^u_{b'}(t)
			&=\PP\left(
			\tau_{\Pi^+(u,b)} \leq t \right)
			-\PP\left(
			\tau_{\Pi^+(u,b')} \leq t 
			\right)\\
			&=\PP\left(\phi_t^u > b \right)
			-\PP\left(\phi_t^u \geq b'\right)\\
			&= \PP\left(\phi_t^u \in (b,b')\right)\\
			&\geq \PP\left(\phi_t^u \in (b,b+\epsilon)\right)\\
			&>0,
		\end{align*}
		which implies $ F^u_b(t)>F^u_{b'}(t) $.
		
		Consider then some proper contour $ \partial\B $, and assume, for contradiction, that we have $ B(\B,u) > C_Q(u) $ for some $ u\in S^{d-1} $. Since, by definition, $ C_Q(u)\in \U^u_{t_s} $ we must also have $ F^u_{C(u)}(t_s)>F^u_{B(\B,u)}(t_s) $ which yields
		\begin{align*}
			\PP(\tau_{\Pi^+(u,B(\B,u))} > t_s)
			&=  1-F^u_{B(\B,u)}(t_s)\\
			&>  1-F^u_{C_Q(u)}(t_s)\\
			&= \PP(\tau_{\Pi^+(u,C_Q(u))} > t_s)\\
			&= q_s.
		\end{align*}
		This contradicts the fact that $ \partial\B $ is a proper contour and we must therefore have $ B(\B,u) = C_Q(u) $ for all $ u\in S^{d-1} $.
	\end{proof}
\end{prop} 

We can also extend this result to proper contours in the return period sense.

\begin{prop}
	Fix some $ t_r $ such that $ C_T(u) $ is defined for all $ u\in S^{d-1} $ and let the conditions of \Cref{prop:propQ} hold. We also assume that $ \tau_{\Pi^+(u,C(u))} $ is non-deterministic in the sense that $ \PP(\tau_{\Pi^+(u,C(u))}=t_r)<1 $.
	
	Under these conditions, if there exists some proper contour $ \partial\B $ in the return period sense, i.e.\ $ \E\left[\tau_{\Pi^+(u,B(\B,u))}\right]=t_r $ for all $ u\in S^{d-1} $, then $ B(\B,\cdot)=C_T $ and
	$$\B=\bigcap_{u\in S^{d-1}}\Pi^-(u,C_T(u)).$$
	\begin{proof}
		Consider the proper contour $ \partial\B $ and assume for contradiction that $ B(\B,u) > C_T(u) $ for some $ u\in S^{d-1} $. We then have that since $ \TT_u(B(\B,u)) = \TT_u(C_T(u)) $ then
		\begin{align*}
		0
		&=\TT_u\big(C(u)\big)-\TT_u\big(B(\B,u)\big)\\
		&=  \int_{0}^{\infty}\left(1-F^u_{C_T(u)}(t)\right)dt
		 -  \int_{0}^{\infty}\left(1-F^u_{B(\B,u)}(t)\right)dt\\
    	&=  \int_{0}^{\infty}\left( F^u_{B(\B,u)}(t)-F^u_{C_T(u)}(t)  \right)dt.
		\end{align*}
		Since $ F^u_{B(\B,u)}(t)-F^u_{C_T(u)}(t) \geq 0 $ this equality implies $ F^u_{B(\B,u)}(t)=F^u_{C_T(u)}(t) $ for almost all $ t\in[0,\infty) $. However, $ b\mapsto F^u_b(t) $ is monotone decreasing on $ \U^u_{t} $, so we must have $ C_T(u)\notin\U^u_{t} $, implying $ F^u_{C_T(u)}(t)\in \{0,1\} $, for almost all $ t\in[0,\infty) $. Furthermore, since $ t\mapsto F^u_b(t) $ is monotone non-decreasing $ F^u_{C_T(u)}(t) $ either equals $ \1(t<s)$ or $ \1(t\leq s)$ for some $ s\in[0,\infty) $, where $ \1 $ denotes the indicator function. In fact, since $ \int_{0}^{\infty}\left(1-F^u_{C_T(u)}(t)\right)dt=t_r $ we get $ F^u_{C_T(u)}(t)= \1(t<t_r) \text{ or } \1(t\leq t_r)$. From this we see that $ \PP(\tau_{\Pi^+(u,C(u))}=t_r)=1 $ which contradicts our assumption that $ \tau_{\Pi^+(u,C(u))} $ is non-deterministic, thereby implying $ B(\B,u) = C_T(u) $ for all $ u\in S^{d-1} $.
	\end{proof}
\end{prop}

These results shows us that any proper contour, in either the return period or the quantile sense, is uniquely defined by \eqref{eq:propconstr} with $C=C_Q$ or $C_T$ respectively.

In the case where no proper contour exists we will need alternative methods for constructing a valid contour. One example of a possible construction is the following. Here $x\in\R^d$ is some suitable centre point, ideally such that $C(u)-\langle u,x\rangle > 0 $ for all $u\in S^{d-1}$, around which the contour is drawn.

\begin{equation}\label{def:convhB}
\B= \clos \left(\convh\left(\left\{x+u\left(C(u)-\langle u,x\rangle\right)^+:u\in S^{d-1}\right\}\right)\right) ,
\end{equation}
where $(\cdot)^+$ equals $\max(\cdot,0)$, $ \convh(\cdot) $ denotes the convex hull and $\clos(\cdot)$ is the closure.

\begin{prop}\label{prop:validcont}
	Fix some $ t_s \in (0,\infty),\,q_s\in (0,1) $ such that $ C_Q(u) $ is well defined for all $ u\in S^{d-1} $ and bounded from above. Let $ \B $ be constructed as in  \eqref{def:convhB} with $C=C_Q$, we then have that $ \partial\B $ is a valid contour in the quantile sense.
	
	Similarly, if $C_T$ is defined and bounded above for some $t_r\in (0,\infty)$ and $\widehat{\B}$ is constructed as in \eqref{def:convhB} with $C=C_T$, we then have that $ \partial\widehat{\B} $ is a valid contour in the return period sense.
	\begin{proof}
		Since, by definition, $ x+u\left(C_Q(u)-\langle u,x\rangle\right)^+\in\B $ we have
		$$  B(\B,u) = \sup\{\langle u,v\rangle :v\in\B\} 
		\geq \big\langle u,x+u\left(C_Q(u)-\langle u,x\rangle\right) \big\rangle 
		=C_Q(u),  $$
		for any $ u\in S^{d-1} $. Lastly, $ \B $ is closed and convex by definition, and since $C_Q$ is bounded from above we have that $ \left(C_Q(u)-\langle u,x\rangle\right)^+$ is bounded. As a consequence, $ \B $ is compact, making $\partial\B$ valid in the quantile sense. An identical argument shows that $\partial\widehat{\B}$ is valid in the return period sense.
	\end{proof}
\end{prop}

To ensure that this construction is always feasible we have the following result, which ensures that $C_Q$ and $C_T$ are indeed bounded, thereby guaranteeing the existence of valid contours.

\begin{lemm}\label{lem:Cbounded}

	Fix some $ t_r,\,t_s\in (0,\infty),\, q_s\in(0,1) $ such that $C_T$ is well defined and assume that \eqref{axiom:contmax} holds. We then have that $ C_Q $ and $C_T$ are bounded from above on $ S^{d-1} $.

	\begin{proof}
			Firstly, by \Cref{thm:Fbtcont}, we have that $ C_Q $ is defined and finite for any $ u\in S^{d-1} $. Furthermore, if we define $\phi_t^\infty=\sup_{s\in [0,t]}\|V_s\|$ we may note that \eqref{axiom:contmax} implies $\PP(\phi_t^\infty=\infty)=0$ for any $t\in (0,\infty)$. 
			
			As a consequence we may pick some $b\in\R$ such that $\PP(\phi_{t_s}^\infty<b)>q_s$ and compute $$\PP(\tau_{\Pi^+(u,b)}>t_s)\geq \PP(\phi_{t_s}^\infty<b) > q_s.$$ Due to this, we see that $C_Q(u)<b$ for all $u\in S^{d-1}$.
			
			Similarly, we may pick $b'\in\R$ such that $\PP(\phi_{2t_r}^\infty<b')>0.5$. If we then define $K=\{v\in\R^d:\|v\|\geq b'\}$ we get 
			$$ 
			\E\left[\tau_{\Pi^+(u,b')}\right]\geq
			\E\left[\tau_{K}\right]\geq
			2t_r\PP(\tau_{K}>2t_r)\geq
			2t_r\PP(\phi_{2t_r}^\infty<b')>
			t_r.
            $$
			This implies $C_T(u)<b'$ for all $u\in S^{d-1}$.

		\end{proof}
\end{lemm}

With this result we can guarantee the existence of valid contours. However, other construction methods also exist. In \cite{convcont}, the authors consider a scenario where $C_e$ could produce a proper contour, but due to estimation errors, the approximated $C_e$ fails to do so. To address this issue, they propose constructing an inflated contour using \eqref{eq:propconstr} based on $C_e+c$ for some appropriate $c\in\R$. In a more general case where $C_e$ does not admit a proper contour, as presented in \cite{voronoi}, an invalid contour is constructed using \eqref{eq:propconstr}, followed by an extension procedure that guarantees a valid construction in the limit. For both cases it is assumed that $C_e$ is bounded, and hence \Cref{lem:Cbounded} ensures that these methods can still be applied in our setting. We also mention \cite{sande2023minimal}, where a numerical algorithm for computing a minimal valid contour in the sense of mean width is developed and analysed.

With the existence and construction of contours settled we can move on to some examples. The goal of these are to show the ways the presented methods differ from the existing framework presented in e.g.\ \cite{firstaltcontour}. The first main difference is the ability to consider a continuous-time framework which may more accurately capture the dynamics of $ V $ in between the discrete points. Secondly we can also allow non-stationary behaviour which allows the inclusion of effects like climate change as a part of the model for $ V $.

\section{Theoretical Example} \label{sec:thex}

We will in this section aim to define a proper contour in the return period sense with a target return period of $ t_r $, under the assumption that $V$ follows the continuous dynamics described below. Once this is done we will compare our exact method with an i.i.d.\ method to highlight the differences. In analysing these models we also discuss when i.i.d.\ methods produce conservative estimates.

It is also worth mentioning that a similar case was studied in \cite{Leira}. Here, several types of contours were compared for both continuous and discrete models of the underlying environmental factors. It was also observed that a continuous model for $V$ may result in larger contours. Additionally, in \cite{huseby2023AR1,mackay2021effect,vanem2023analysing} they compare discrete models taking serial correlation into account with i.i.d.\ models. In these articles they conclude that including autodependence will lead to smaller contours, and hence argued that ignoring correlation will lead to overly conservative estimates. {The i.i.d.\ method will therefore serve as a conservative representation of all possible discrete models. We can therefore compare its resulting contours to that of a continuous model to examine how autodependence affects the size of the contour in a continuous-time setting.}

We here consider the case where $V_t \in \R^1$ is defined by 
$$V_t = \sqrt{2\theta} \int_{-\infty}^{t}e^{-\theta(t-s)} dW_s,$$
where $ W $ is standard Brownian motion and $ \theta\in\R,\,\theta>0 $. This makes $ V $ a standardised Ornstein-Uhlenbeck  process which serves as a continuous interpolation of an AR(1) discrete-time process. Note that $V$ is standardised to ensure a mean of 0 and variance 1, which implies $V_t$ standard normally distributed for any $t$. {This allows us to interpret the continuous-time process in the sense of \Cref{sec:response}. We assume that the long-term conditions over the periods $[n\Delta t,(n+1)\Delta t]$, $n\in\N$, are standard normally distributed and follow an AR(1) process. Consequently, $V_t$ is an interpolation of these conditions and can be interpreted as the average conditions over $[t,t+\Delta t]$.}

Since $ -V $ satisfies $ -V_t = \sqrt{2\theta} \int_{-\infty}^{t}e^{-\theta(t-s)} d(-W)_s$, we see that $ -V $ is also an Ornstein-Uhlenbeck process with the same parameters and thus equal in law to $ V $. As such we know that $ C_T $ is constant on $ S^0=\{-1,1\} $. In fact, if we considered a $ d $-dimensional Ornstein-Uhlenbeck process, $ C_T $ would still be constant and $ \partial\B $ would equal a $ (d-1) $-sphere with a radius given by the same value of $C_T$ as our $1$-dimensional case.

In computing $ C_T $ we again consider $ \TT_u(b)= \E\left[\tau_{\Pi^+(u,b)}\right] $. Under the assumption that $ V_0=0 $ we have an explicit representation of $ \TT_u(\cdot) $, independent of $u$, given in e.g.\ \cite{OUmom} as
$$ \TT_u(b)= \frac{1}{2\theta}\sum_{i=1}^{\infty} \frac{ (\sqrt{2}b)^{i}}{i!}\Gamma\left(\frac{i}{2}\right), $$
where $ \Gamma $ is the gamma function. One can also show the more convenient alternative representation of
$$ \TT_u(b)= {\frac{\sqrt{\pi}}{\theta\sqrt{2}}}\int_{0}^{b}\left(1+\text{erf}\left(\frac{t}{\sqrt{2}}\right)\right)e^{t^2/2}dt, $$
where erf denotes the error function. By inverting $ \TT_u(\cdot) $ we can easily compute $ C_T $ numerically. The resulting contour would then, due to $ C_T $ being constant, be the two points given by $ \partial\B=\{\pm \TT_u^{-1}(t_r)\} $.

This gives us an explicit representation of the optimal contour. However, as an alternative we could consider am i.i.d\ model for $ V $. We define $\overline W= \{\overline W_n\}_{n=0}^\infty $ as an independent sequence of standard normally distributed random variables. Further define $ \overline{V} $ by $ \overline{V}_t=\overline W_{\lfloor t/\Delta t \rfloor} $, for some $ \Delta t >0 $ and note that $ V_t $ and $ \overline{V}_t $ are equal in law for any $ t\in\R $. {As such we may consider $ \overline{V} $ as an alternative i.i.d.\ model of $ V $.}

If we were to apply the i.i.d.\ method presented in e.g.\ \cite{firstaltcontour} we could compute $ C_T $ based on $ \overline{V} $. We know that this is equivalent to considering $C_e$ with an exceedence probability of $ p_e=\Delta t /t_r $, which means that for all $ u\in S^{0} $ we would have
$$  C_e(u) =  \Phi^{-1}\left(1-\frac{\Delta t}{t_r}\right), $$
where $ \Phi $ is the cumulative distribution function of a standard normal random variable.

With these two models we can compare how the resulting contours differ. Since both models provide perfectly circular contours (insofar as $S^0$ can be referred to as \textit{circular}) we can instead compare the radii. To compute exact numbers we want a reasonably realistic value for $ \theta $. This will be chosen based on a time series, $\{H_n\}_{n\in\N}$, of significant wave heights. The details of this dataset will be given in \cref{sec:empex}.

{
In order to choose our parameters we first pick $ \Delta t = 3 $ hours and compute the $24$ hour autocorrelation ($ AC_{24} \approx 68\%$) of the standardised data $ (H-\mu_H) /\sigma_H $. Here $ \mu_H $ and $ \sigma_H $ are the empirical mean and standard deviation, respectively, of the time series $ \{H_n\}_{n\in\N} $. By noting that the $24$-hour autocorrelation of an Ornstein-Uhlenbeck process satisfies $ \theta=-log(AC_{24})/24 $ we can estimate $ \theta\approx 0.16 $ hours$^{-1}$.


} 

\begin{figure}[h]
	\includegraphics[width=0.95\linewidth]{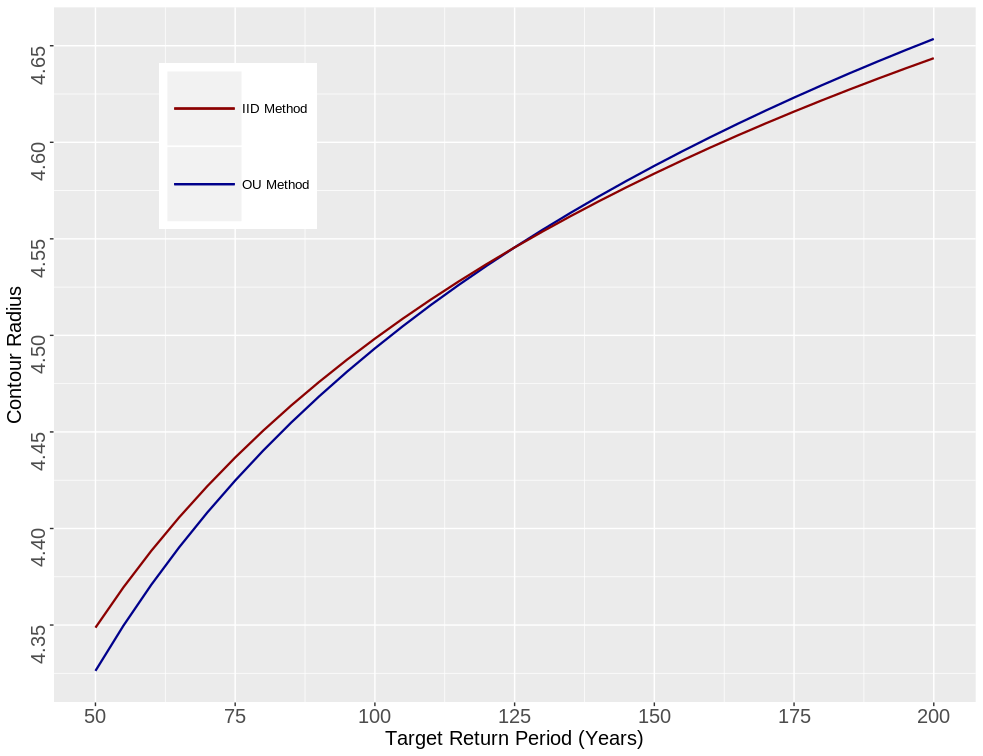}
	\caption{Comparison of contour radii for different methods}
	\label{fig:OU_comparison}
\end{figure}

The resulting radius curve based on the continuous-time model of $ V $, labelled \textit{OU Method}, and the curve based on $ \overline{V} $, labelled \textit{IID Method}, are plotted in  \Cref{fig:OU_comparison}.

\begin{rem}
    The differences are quite small, but there is still a noticeable distinction between the two methods. { In fact, their similarities are heavily dependent on the value of $\theta$,  a higher value would push the radius of the O.U.\ method more noticeably above or below the one of the i.i.d.\ method. For example, by accounting for trends and seasonality in $H$ one may get $\theta\approx 0.025$. This would yield  e.g.\ a 200-year contour radius of 4.75, which is a 2.3\% increase over the i.i.d.\ method. If we instead considered an estimate of $\theta$ based on e.g.\ the 7-hour empirical autocorrelation, then $\theta\approx 0.01$ is obtained. This estimate would yield a 200-year contour radius of 4.54, a 2.3\% decrease compared to the i.i.d.\ method.}     
\end{rem}

As we see the i.i.d.\ method produces larger and more conservative contours for low values of $ t_s $, but crosses below for sufficiently large return periods. One can even find an approximation of when the two lines cross by noting that this point occurs when the radius, $ R $, of the contour satisfies
\begin{equation}\label{eq:crossR}
	\frac{\Delta t}{1-\Phi(R)}
	=
	\TT_u(R).
\end{equation}
It can further be shown, by taking specific asymptotic expansions, that for high values of $ R $ we have the following approximations
$$ 
\frac{\Delta t}{1-\Phi(R)} \approx \Delta t \sqrt{2\pi}Re^{R^2/2},
\quad
\TT_u(R) \approx \frac{\sqrt{2\pi} e^{R^2/2}}{\theta R}.
$$
Applying these to both sides of \eqref{eq:crossR} then yields
$$
\Delta t \sqrt{2\pi}Re^{R^2/2}
\approx
\frac{\sqrt{2\pi} e^{R^2/2}}{\theta R},
$$
as long as the point where the lines cross occur for sufficiently high values of $ R $. Simplifying this expression we get the approximate identity $ \theta\Delta t R^2=1 $, which means we can compute the return period for which this point occurs, here denoted $ t_r^* $, by
\begin{equation}\label{eq:crosstime}
	t_r^*
	\approx
	{\sqrt{\frac{2\pi\Delta t}{ \theta}} e^\frac{1}{2\theta\Delta t}}.
\end{equation}
For our specific parameters, this approximation yields $ t_r^*\approx 130$ years. 

By analysing \eqref{eq:crosstime} we can supplement the more heuristic reasons why these different models yield different contours.

Using $ \overline{V} $ ignores the autocorrelation, allowing the process to vary more wildly, usually producing larger contours. This effect is magnified when $ \Delta t $ is low which corresponds to the limit $ \lim_{\Delta t\to 0} t_r^* = \infty $.  This effect justifies why the i.i.d.\ approach is usually considered as a conservative estimate, as argued in e.g.\ \cite{huseby2023AR1,mackay2021effect,vanem2023analysing}. In these articles they remark that including serial correlation would lead to less conservative contours, and thus less extreme design conditions. This holds in the setting of discrete models for $V$. However, as we see from this example, there are situations where the consideration of continuous-time models can lead to more conservative contours despite the inclusion of serial correlation.

{
As this example shows, there are situations where the i.i.d.\ method can overestimate the return period compared to a continuous-time model. Since the presumed true model of $ V $ is non-discrete, it has the possibility of exceeding boundaries at times between the discrete points. Since $V_t$ was calibrated to equal $\overline{V}_t$ in distribution, and therefore represents the long-term conditions over $[t,t+\Delta t]$, we consider the following scenario. Some off-shore structure is exposed to environmental loads, for simplicity we focus only on the significant wave height. The significant wave heights over the periods $[0,\Delta t]$ and $[\Delta t,2\Delta t]$, each induces their own short-term probability of failure. However, due to the variability of $V$, the significant wave height could potentially be higher over the middle period of $[\Delta t/2,3\Delta t/2]$, than over any of the two other periods individually. This permits a potentially much higher failure probability over the middle period that would be ignored by only considering $[0,\Delta t]$ and $[\Delta t,2\Delta t]$. As a consequence, the time until failure could be underestimated by not considering the short-term variability of the long-term conditions described by $V$. This effect is reflected by $ \partial t_r^*/\partial\theta < 0 $, which implies that increasing $ \theta $ shrinks the domain where the method based on $ \overline{V} $ is conservative. Indeed, a high $ \theta $ increases the short-term volatility of $ V $, thus improving the chances of the process exceeding a fixed boundary at times in-between the points of $ \{n\Delta t, n\in\N\} $.
}

{For the purposes of response analysis we usually want to use a discrete model for $V$. This example, however, shows that this can overestimate the return period of sufficiently extreme conditions. Consequently, those models may underestimate the response resulting from such conditions. 
Bear in mind that the i.i.d.\ method can be viewed as a conservative representative for all possible discrete methods. 
Thus, if the continuous-time modelling of $V$ produces meaningfully larger contours than even the conservative i.i.d.\ estimate,  then this serves as an indication that a lower $\Delta t$ should be considered.}

\section{Empirical Example} \label{sec:empex}

\subsection{Data and Outline of Method}

The data considered for this example will consist of ERA5 reanalysis data \cite{ERA5}. We will use hourly data for significant wave height and wave periods from 65\degree N, 0\degree W over the period $ 1959$-$2021 $. While the calibration will use the full resolution of one hour, we will use a three-hour time step for the purposes of simulation and computation.

Our primary goal is here to compute an empirical environmental contour in the quantile sense for a survival time of $ t_s= 50 $ years and a survival probability of $ q_s = e^{-1} \approx 37\% $. In doing so we will also present a specific algorithm for generating such contours which is carried out in four steps:
\begin{itemize}
	\item The distribution of $V_t$ is estimated for all relevant values of $t$.
	\item Paths of $V$ are simulated by using sequences of independent (but not identically distributed) random variables.
	\item These paths are used to compute samples of $\phi_{t_s}^u$ which allows the computation of $C_Q$ by \Cref{lem:phieq}.
    \item The contour can then be constructed. If a proper contour exists, then it is defined by \eqref{eq:propconstr}, but if no such contour exists, we may employ the methods of \cite{sande2023minimal} to produce a minimal valid contour in the sense of mean width.
\end{itemize} 

Note that a survival probability of $e^{-1}$ is chosen to correspond with a 50-year return period if $\tau$ had an exponential density. Due to the presence of a trend we do not have this distribution, but it will still serve as an easy, although rough, point of comparison.

\begin{rem}
    This procedure is very similar to the method proposed in \cite{vanem2023analysing}, which consider a discrete, but serially correlated model for $V$. The author also considers contours in the return period sense in the same manner as defined in \Cref{sec:dyncontour}. Firstly, both the marginal distribution and autocorrelation function of $V$ are estimated. It is then assumed that $V_t=g(Z_t)$, where $g$ is the Rosenblatt transform, implying that $Z$ has standard Gaussian marginals. The author then finds a unique autocorrelation structure for $Z$ that implies the empirical autocorrelation of $V$. Paths of $V$ can then be simulated by simulating $Z$, which is used to estimate $C_Q$.

    In our simulation method we consider a non-stationary distribution of $V$, but ignore its autocorrelation. It is entirely possible to include serial correlation by extending the method in \cite{vanem2023analysing} to a non-stationary setting. One could for example consider $V_t=g_t(Z_t)$, where $g_t(\cdot)$ is the Rosenblatt transform corresponding to the marginal density of $V_t$.
\end{rem}

\subsection{Calibration of Distributions}

We here have $ d=2 $ with $ V_t=(P_t,H_t) $ where $ P $ and $ H $ denotes the {wave period} and significant wave height respectively. Following e.g.\ \cite{convcont,firstaltcontour,vanemTrend}, we model $ H $ using a 3-parameter Weibull distribution, and $ P $ with a conditional log-normal distribution. However, due to the presence of long-term trends discussed in e.g.\ \cite{Kushnir,vanemTrend,vanemBayes}, we will model $ V $ as non-stationary to take both this trend, as well as seasonality, into account.

We assume that $ H_t\sim W(\lambda_t,k_t,\theta) $, i.e.\ a $ 3 $-parameter Weibull distribution with scale $ \lambda_t $, shape $ k_t $, and location $ \theta $. Here $ \lambda $ is taken to be on the form $ \lambda_t=(c_1+c_2t)l_t $, furthermore, $ l $ and $ k $ are assumed periodic with a period of one year.

As for the wave period, we assume that $ \left(\log(P_t)|H_t=h\right)  \sim\NN(\mu(t,h),\sigma^2(t,h)) $, i.e.\ a conditional normal distribution with mean $ \mu(t,h) $ and variance $ \sigma^2(t,h) $. Here $ \mu $ and $ \sigma $ are assumed to be on the form $ \mu(t,h)=m(t)+f_\mu(h) $ and $ \sigma(t,h)=s(t)f_\sigma(h) $ where $ m $ and $ s $ are periodic with a period of one year.

{
\begin{rem}
    Several simplifications could potentially be considered here. For example, one may ignore seasonality and use the popular parametric estimates $ f_\mu(h) =a_1 + a_2 h^{a_3} $ and $f_\sigma(h)=b_1 + b_2 e^{b_3h} $ for some constants $a_i,\, b_i,\, i=1,2,3$. This would allow for a simpler and fully parametric model for the distribution of $V_t$.
\end{rem}
}

In order to perform our calibration procedure we will first remark that if $ H_t \sim W(\lambda_t,k_t,\theta)$, we then have for any $ \lambda_t',k_t'\in\R $ that
	
	$$ \frac{(H_t - \theta)}{\lambda'_t} 
	\sim W\left(\frac{\lambda_t}{\lambda_t'},{k_t},0\right)
	,\quad
	{(H_t - \theta)^{k'_t}}
	\sim W\left(\lambda_t^{k'_t},\frac{k_t}{k'_t},0\right).$$

To estimate $ (\lambda_t,k_t,\theta) $ for $ H $ we then do the following.

\begin{itemize}
	\item $ \theta $ is estimated by the minimal measured value (rounded down to 2 significant digits to avoid numerical issues).
	\item The linear trend parameters, $ (c_1,c_2) $, are estimated by linear regression on $ H-\theta $.
	\item $ k_t $ is estimated by inverting the equality
	\begin{equation*}
		\frac{\Gamma(1+1/k_t)^2}{\Gamma(1+2/k_t)}= \frac{\E[H_t-\theta]^2}{\E[(H_t-\theta)^{2}]},
	\end{equation*}
	where the expectations are computed by smoothing spline regression of $ H-\theta $. We then normalise $ k $ by defining $ k_t'=k_t/\int_{0}^{1}k_t dt $ so the average value of $ k' $ equals $ 1 $. This is done to avoid numerical issues from taking high powers.
	\item For calibration of $ l $ we note that
	$$
	\E\left[\frac{(H_t - \theta_t)^{k'_t}}{(c_1+c_2t)^{k_t'}}\right]
	=
	l_t^{k'_t}\Gamma\left(1+\frac{k_t}{k'_t}\right).
	$$
	We can then fit $ l_t^{k'_t} $ by a smoothing spline regression of 
	$$ \frac{(H_t - \theta_t)^{k'_t}}{(c_1+c_2t)^{k_t'}\Gamma(1+\frac{k'_t}{k_t})}. $$
\end{itemize}

\begin{figure}[h]
	\centering
	\begin{subfigure}{.5\textwidth}
		\centering
		\includegraphics[width=.95\linewidth]{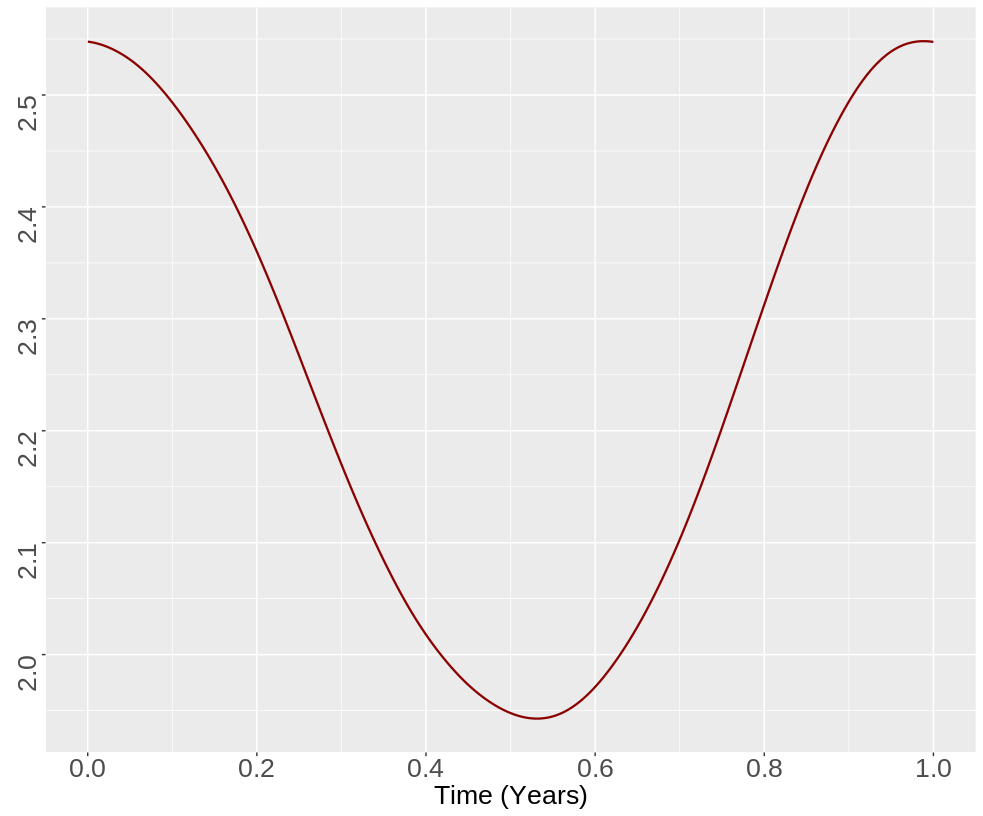}
		\caption{Values of $ k(t) $}
		\label{fig:ks}
	\end{subfigure}%
	\begin{subfigure}{.5\textwidth}
		\centering
		\includegraphics[width=.95\linewidth]{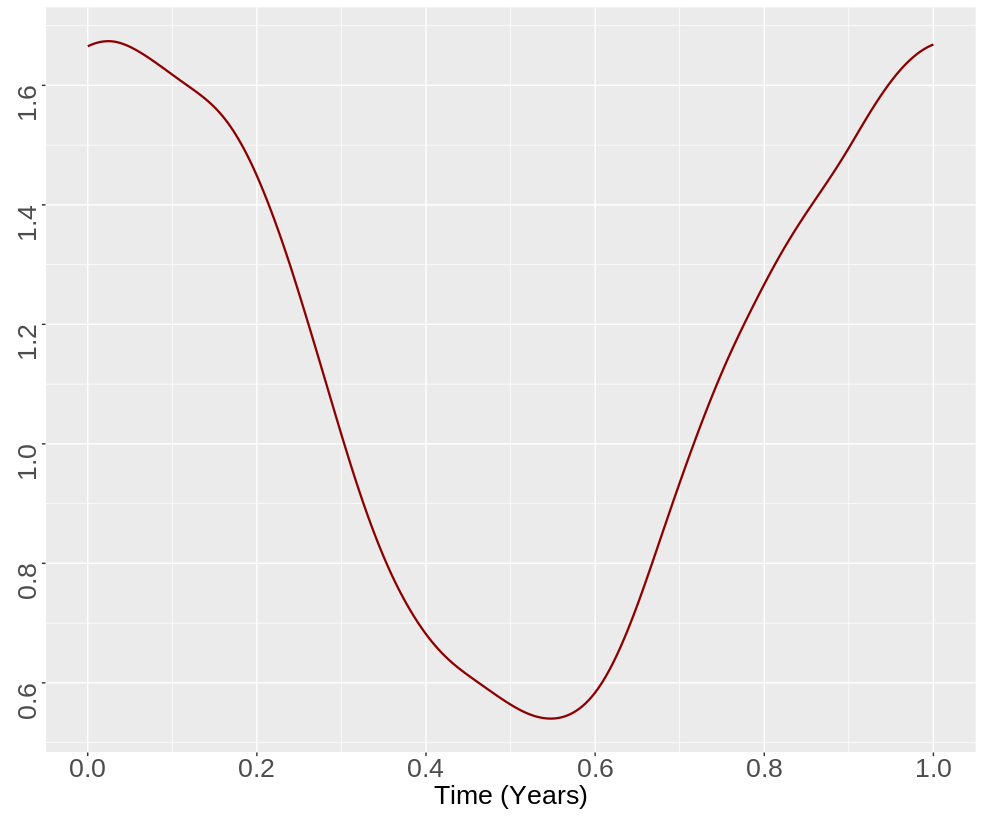}
		\caption{Values of $ l(t) $}
		\label{fig:lams}
	\end{subfigure}
	\caption{Non-parametric functions for $ H $}
	\label{fig:Hfuncs}
\end{figure}

\begin{table}[h]
	\begin{center}
		\begin{tabular}{|c|c|c|}
			\hline
			$ \theta $ & $ c_1 $ & $ c_2 $\\
			\hline
			 0.37 m  &  2.5 m &  4.0e-3  m/y \\
			\hline
		\end{tabular}
	\end{center}
	\caption{Parameters for $ H $}
	\label{fig:Hpars}	
\end{table}

The resulting parameters are given in \Cref{fig:Hpars} and \Cref{fig:Hfuncs}. Note that the parameters in \Cref{fig:Hpars} are given in meters (m) and meters per year (m/y), the functions in \Cref{fig:Hfuncs} are dimensionless. Furthermore, we consider $ t=0 $ to occur at the end of the dataset, i.e.\ the beginning of 2022.

As for $ P $ we do the following:
\begin{itemize}
	\item We first estimate $ \mu(t,h) $ by $\E\left[\log(P_t)|H_t\right]= \mu(t,H_t)=m(t)+f_\mu(H_t)$. With this we can fit smoothing splines for $ m $ and $ f_\mu $ by generalised additive model calibration.
	\item Similarly, $ \sigma(t,h)=s(t)f_\sigma(h) $ can be computed by 
	$$\E\left[\log\left(\log(P_t)-\mu(t,H_t)\right)^2\right]=L+\log \left(s^2(t)\right)+\log \left(f^2_\sigma(H_t)\right), $$
	where $ L $ is the log-moment of a chi-squared random variable with $ 1 $ degree of freedom. This allows us to fit smoothing splines for $ \log(f^2_\sigma(h)) $ and $ \log(s^2(t)) $ by a weighted generalised additive model calibration. Finally, to counteract issues arising from log-scale calibration, $ s $ is scaled to ensure that $ \left(\log(P_t)-\mu(t,H_t)\right)/\sigma(t,H_t) $ has a variance of $ 1 $.
\end{itemize}

\begin{figure}[h]
	\centering
	\begin{subfigure}{.5\textwidth}
		\centering
		\includegraphics[width=.95\linewidth]{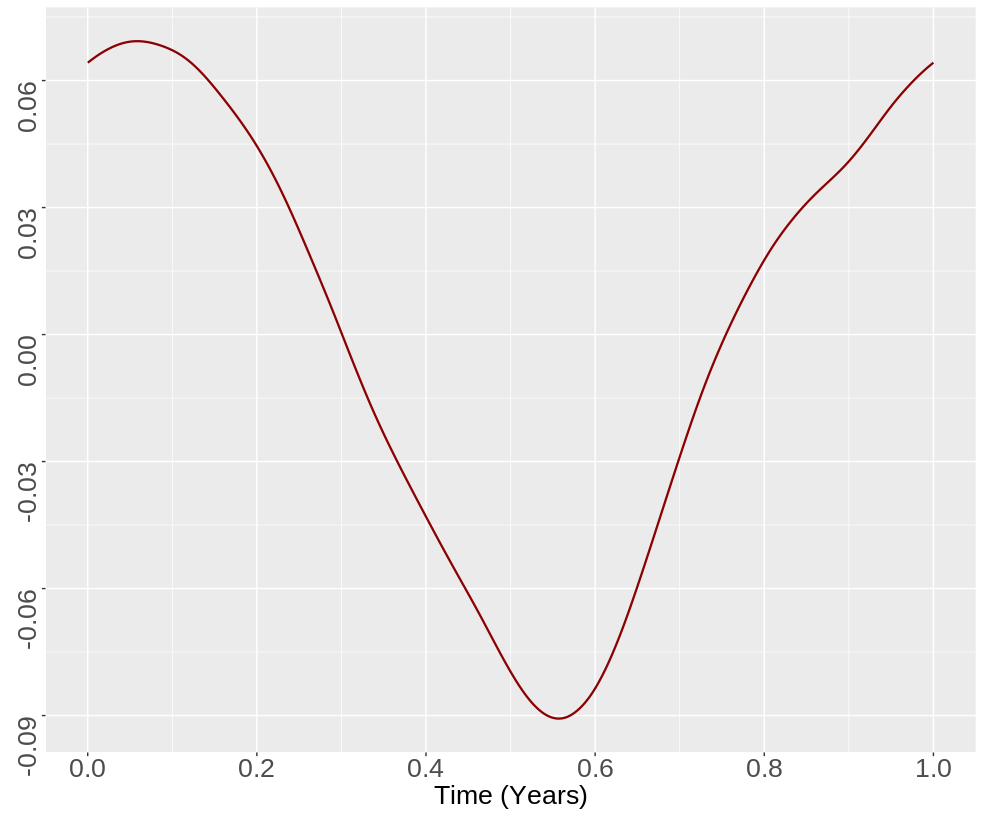}
		\caption{Values of $ m(t) $}
		\label{fig:mus}
	\end{subfigure}%
	\begin{subfigure}{.5\textwidth}
		\centering
		\includegraphics[width=.95\linewidth]{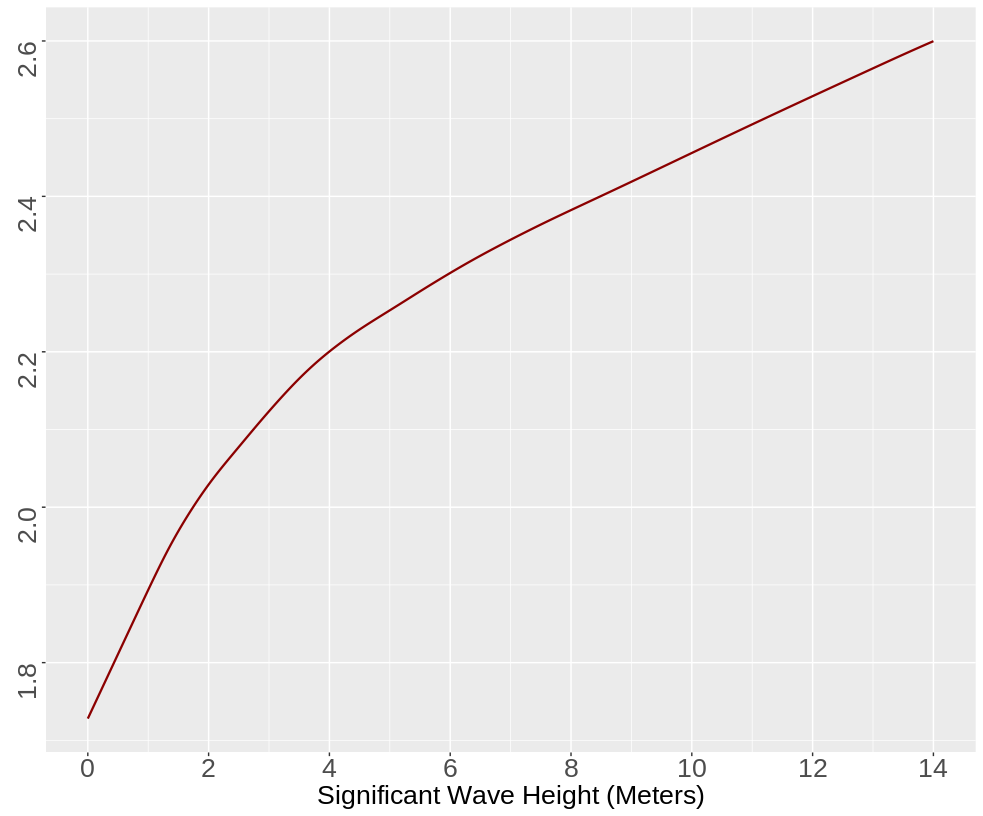}
		\caption{Values of $ f_\mu(h) $}
		\label{fig:muh}
	\end{subfigure}%
	
	\begin{subfigure}{.5\textwidth}
	\centering
	\includegraphics[width=.95\linewidth]{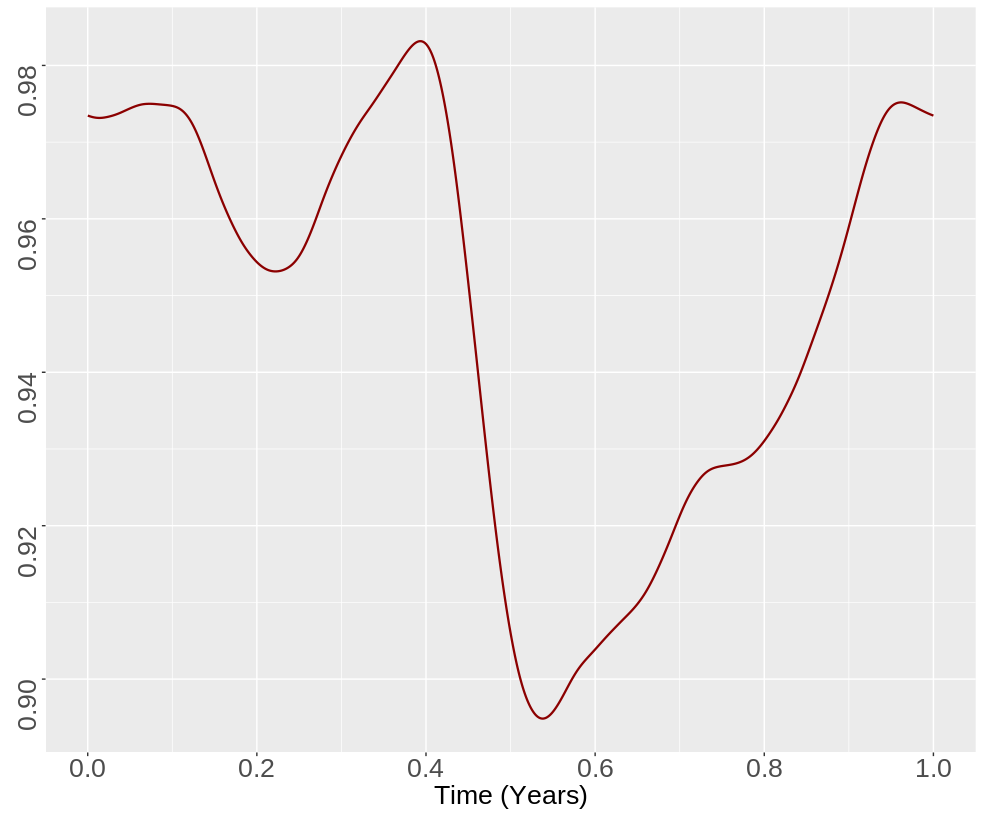}
	\caption{Values of $ s(t) $}
	\label{fig:sigs}
	\end{subfigure}%
	\begin{subfigure}{.5\textwidth}
	\centering
	\includegraphics[width=.95\linewidth]{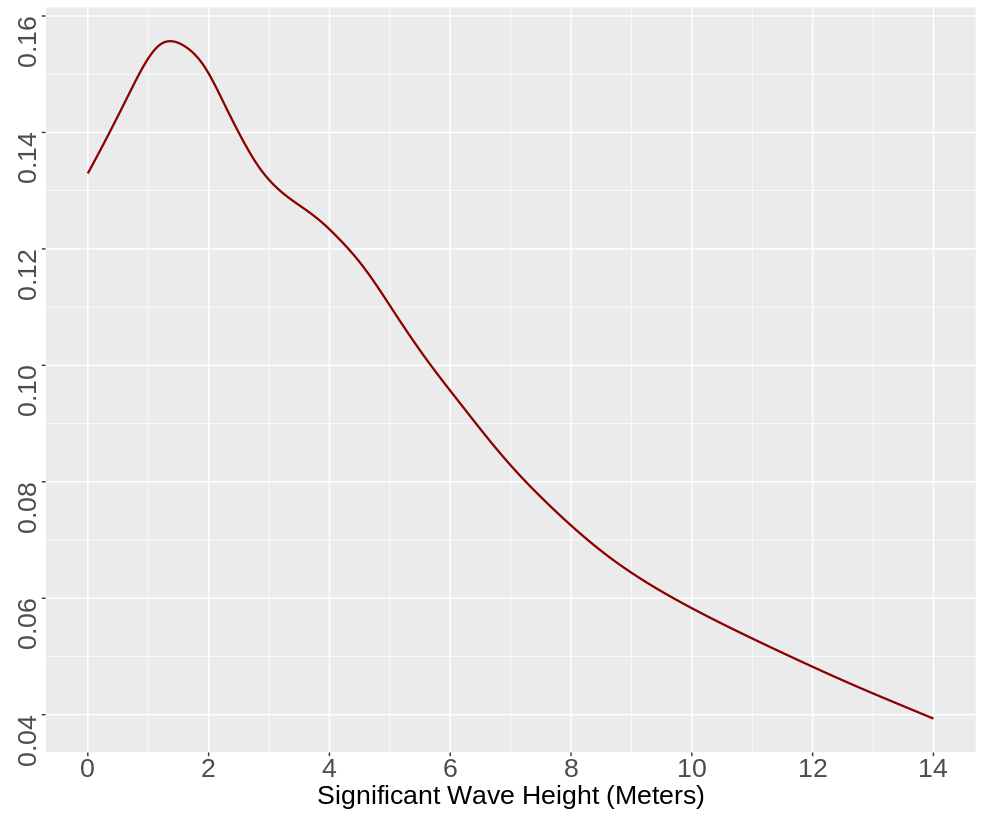}
	\caption{Values of $ f_\sigma(h) $}
	\label{fig:sigh}
	\end{subfigure}
	\caption{Non-parametric functions for $ P $}
	\label{fig:Pfuncs}
\end{figure}

The resulting functions are given in \Cref{fig:Pfuncs}.

\subsection{Simulation}
 
With the $d$-dimensional marginal distributions of $ V $ determined we can move on to the computation of $ C_Q $. Similarly to the previous example we define $ \{W_n\}_{n=0}^\infty $ as a sequence of independent random variables such that $ W_n=V_{n\Delta t} $ in distribution for $ \Delta t=3 $ hours. This again lets us model $ V $ by $ {V_t} = W_{\lfloor t/\Delta t \rfloor} $.

With this model we can easily simulate paths of $ {V} $ over the next $ 50 $ years. Based on these simulations we obtain samples of $\phi^u_{t_s}$ for $180$ uniformly spaced unit vectors in $S^{1}$. By considering the lower $q_s$ quantile of $\phi^u_{t_s}$ for a fixed $u\in S^1$ we obtain $C_Q(u)$ by \Cref{lem:phieq}, which yields $q_s=\PP(\tau_{\Pi^+(u,C_Q(u))}>t_s)=\PP(\phi_{t_s}^u\leq C_Q(u))$.

{
\begin{rem}
    We may also note that $W_n$ has a continuous density for all $n\in\N$. This implies that, for any $u\in S^1, t_s>0$, $\phi^u_{t_s}$ has a density as well. Consequently, by \Cref{thm:Fbtcont}, $C_Q$ is well defined, and by \Cref{prop:validcont,lem:Cbounded} we can guarantee existence of valid contours. Furthermore, as the density of $W_n$ has support equal to $\R\times [\theta,\infty)$, we know that $\phi^u_{t_s}$ has a density with connected support. Consequently, if a proper contour exists, \Cref{prop:propQ} implies it is defined by \eqref{eq:propconstr}.
\end{rem}
}

In existing literature, such as \cite{firstaltcontour} and \cite{vanemTrend}, the inclusion of climatic trends is accomplished by considering a stationary distribution with parameters modified to reflect the observed trend. In order to study the effects of replacing the trend by adjusting parameters of a stationary model, we will examine three cases. This will be done by fixing the trend at either the beginning or end of our 50-year period.
\begin{itemize}
	\item Case $ 1 $: $ \lambda_t = (c_1 + c_2*50y)l_t$, which represents sea-states based on the trend 50 years after 2022.
	\item Case $ 2 $: $ \lambda_t = (c_1 + c_2*t)l_t$, which represents the estimated true sea-state distributions.
	\item Case $ 3 $: $ \lambda_t = c_1l_t$, which represents sea-states based on the trend at the beginning of 2022.
\end{itemize}
Note that we here still include the seasonal effects. This is done to avoid using several calibration methods, which could create artificial differences between the cases unrelated to the trend. However, despite seasonality, we would still expect case 1 and 3 to properly represent stationary alternatives to our method.

\begin{figure}[h]
	\centering
	\begin{subfigure}{.5\textwidth}
		\centering
		\includegraphics[width=.95\linewidth]{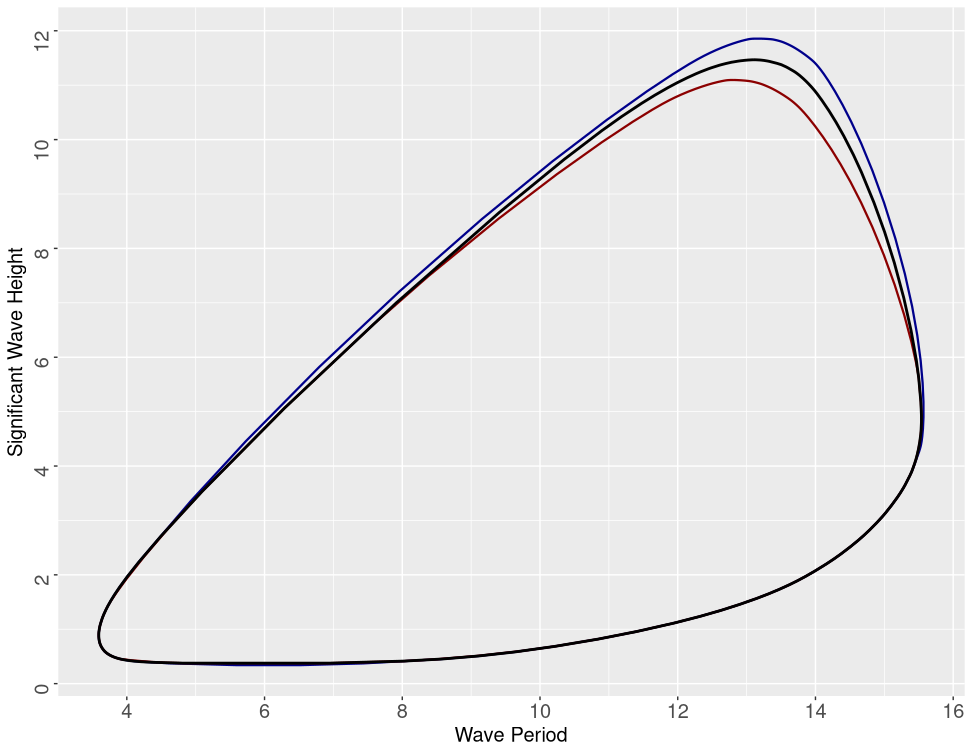}
		\caption{Full Contour}
		\label{fig:Bcasefull}
	\end{subfigure}%
	\begin{subfigure}{.5\textwidth}
		\centering
		\includegraphics[width=.95\linewidth]{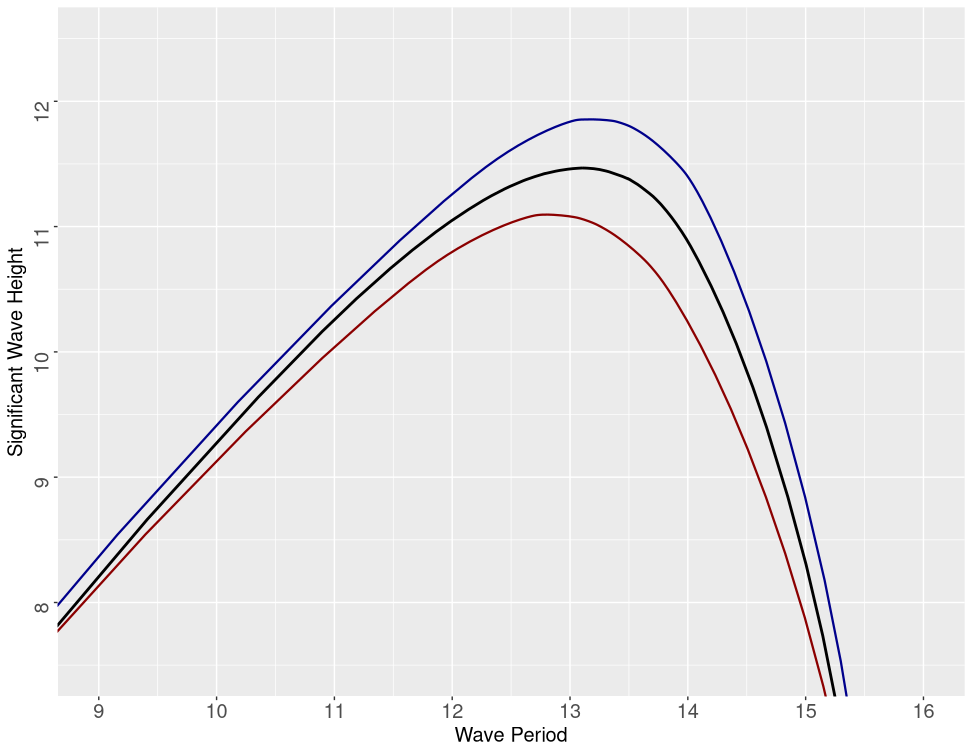}
		\caption{Zoomed-In Contour}
		\label{fig:Bcasezoom}
	\end{subfigure}
	\caption{Environmental contours for case 1 (blue), case 2 (black), and case 3 (red)}
	\label{fig:Bcases}
\end{figure}

As we see, there is a notable difference in $C_Q(u)$, though largely for $u\approx (0,1)$. Note that since $\langle V_t,(0,1)\rangle=H_t$ these values depend mostly on the behaviour of the significant wave height. This demonstrates that including trends is important to avoid underestimation of risk, such as in case 3. However, by including the trend as a non-constant effect, as in case 2, we can still safely reduce the resulting contour relative to case 1, where the highest trend value is applied to the entire period. Specifically, using the conservative estimate of case 1 still overestimates the risk significantly, with a maximal difference in $C_Q$ of 0.41. Since the modelling of $V$ by an i.i.d.\ process is inherently on the safe side, we may be overly cautious by choosing a method which makes further conservative approximations.

\section{Summary and Conclusions}

This paper has rigorously defined and established minimal conditions for existence of environmental contours based on general stochastic processes. These definitions have several advantages over conventional constructions. Chiefly, the ability to properly include climate trends, but also the capability of including seasonality and autodependence.

The theory can also be further generalised to include several extensions. For example, we have discussed the possibility of including serial correlation in the simulation algorithm of \Cref{sec:empex}, by modifying the methods of \cite{vanem2023analysing}. The addition of omission factors can be done in similar way as described in \cite{convcont}. Lastly, buffering, as introduced in \cite{dahl2018buffered}, can be readily extended to our setting, but may require restricting the model choice of $V$, to those based on sequences of independent random variables.

Furthermore, the presented methods have also been compared to conventional techniques, and significant differences in the resulting contours have been demonstrated. In particular, we have discussed how contours can be used to examine the impact of discretisation and autocorrelation. Additionally, we have also illustrated how these methods can avoid the underestimation of risk coming from trends, without the use of excessively conservative strategies. Finally, as part of these examples, we have also presented a strategy for computing these contours based on Monte-Carlo simulation. As such, the approaches considered are presented as an alternative method for the construction of environmental contours.


\section*{Declaration of Competing Interest}
The author declares that there is no known competing financial interest or personal relationship that could have appeared to influence the work reported in this paper
\section*{Acknowledgements}
The author acknowledges financial support by the Research Council of Norway under the SCROLLER project, project number 299897.

The author also want to thank two anonymous referees for providing useful comments and suggestions, which provided significant improvements to the manuscript.\

\bibliographystyle{abbrvdin}
\bibliography{refs}

\end{document}